\documentclass[sn-aps]{sn-jnl}

\usepackage{xspace}
\usepackage{Macros/mydef}
\usepackage{Macros/remarks}
\usepackage{adjustbox}
\usepackage{subcaption}
\usepackage{natbib}

\jyear{2021}%

\theoremstyle{thmstyleone}%
\theoremstyle{thmstyletwo}%

\theoremstyle{thmstylethree}%

\begin{document}

\title[Finite element approximation of MHD]{A high-order residual-based viscosity finite element method for the ideal MHD equations}


\author*[1]{\fnm{Tuan Anh} \sur{Dao}}\email{tuananh.dao@it.uu.se}

\author[1]{\fnm{Murtazo} \sur{Nazarov}}\email{murtazo.nazarov@it.uu.se}

\affil*[1]{\orgdiv{Department of Information Technology}, \orgname{Uppsala University}, \orgaddress{\street{Lägerhyddsvägen 2}, \city{Uppsala}, \postcode{751 05}, \country{Sweden}}}


\abstract{We present a high order, robust, and stable shock-capturing technique for finite element approximations of ideal MHD. The method uses continuous Lagrange polynomials in space and explicit Runge-Kutta schemes in time. The shock-capturing term is based on the residual of MHD which tracks the shock and discontinuity positions, and adds sufficient amount of viscosity to stabilize them. The method is tested up to third order polynomial spaces and an expected fourth-order convergence rate is obtained for smooth problems. Several discontinuous benchmarks such as Orszag-Tang, MHD rotor, Brio-Wu problems are solved in one, two, and three spatial dimensions. Sharp shocks and discontinuity resolutions are obtained.}




\maketitle

\section{Introduction}\label{sec1}

Conservation laws play an important role in modeling and understanding physical processes. For example, conservation of mass, energy, and momentum are fundamental principles in gas and fluid dynamics. For ideal flow, these principles are modeled using the well-known system of Euler equations. Despite a long history and extensive research in understanding and solving Euler equations, many important mathematical and numerical questions remain unsolved. Especially for high-speed flows (higher than the speed of sound), solutions to the Euler system lose their regularity and many complex phenomena such as shocks, strong discontinuities, rarefaction, and contact waves develop. Besides, at very high temperature and energy gas transforms into plasma. Plasma is an ionized gas and the dynamic of plasmas is modeled by adding additional terms and equations involving the magnetic field $\bB$ into the Euler system. This resulting system is usually referred to as Magnetohydrodynamics or MHD in short.

Numerical approximation of MHD started as soon as researchers started to simulate the Euler equations, see \citep{Lax_1954,Brackbill_1976}. Existing numerical methods for compressible Euler equations are hard to extend to approximate the MHD system. If one reason is solving $d$-additional equations, where $d$ is the space dimension, other two main reasons are $(i)$ developing high order robust shock-capturing techniques, $(ii)$ satisfying the divergence-free constraint for the magnetic field, \ie $\DIV \bB = 0$. State-of-the-art numerical methods to solve the MHD equations are finite differences, finite volumes, and discontinuous Galerkin (DG) schemes. Most of the finite difference and finite volume schemes are based on approximate solutions to the Riemann problem of MHD, see for instance \citep{Brio_Wu_1988, Powell_et_al_1991,  Dai_Woodward_1998, Balsara_et_al_1999, Bouchut_et_al_2007, Balsara_2010, Balsara_Dubser_Abgrall_2014} and references therein. Due to similarities of finite volume and DG schemes, an approximate Riemann solver approach was incorporated in the Galerkin variational formulation for the DG schemes, see \eg \citep{Warburton_Karniadakis_1999, Li_Shu_2005, Dumbser_2016}. In general, the exact or approximate Riemann solvers are difficult to compute. Especially, the MHD system is not strictly hyperbolic and has a non-convex flux, and therefore it may admit non-regular waves which may result in non-unique Riemann solution \citep{Torrilhon_2002}. To overcome this difficulty, numerical methods that do not use Riemann solvers have been developed. For instance central DG schemes reported in \citep{Tadmor_mhd1_2004, Li_Xu_Yakovlev_2011, Cheng_et_al_2013} approximate the MHD solutions without using any approximate or exact Riemann solvers.

The success of finite volume and DG schemes has not been observed in a continuous finite element (FE) framework. A reason for it could be the lack of high order, robust explicit stabilization techniques to stabilize continuous FE approximations of MHD. Traditional stabilization methods for FE approximations of conservation laws including Euler and MHD equations are based on the Galerkin-Least-Squares argument (GLS) augmented with residual-based shock-capturing terms, see for instance \citep{Szepessy_1989, Johnson_et_al_1990} and references therein. GLS schemes are implicit by construction, difficult to make high order in time, and nontrivial to implement since the test functions are required to be time-dependent. Only a few works were done to simulate resistive MHD using GLS stabilization, see for instance \citep{Badia_Codina_2013, Sitaraman_Raja_2013, Shadid_et_al_2016} and references therein.

It appears that a way of constructing a high-order explicit FE approximation of conservation laws could be by entirely disregarding the GLS terms from the stabilized formulation. Recently, \cite{Guermond_2008, Guermond_pasquetti_popov_JCP_2011, GuermondNaPo11, Nazarov_Larcher_2017} proposed the so-called entropy viscosity method, where the FE approximation is stabilized by adding an elliptic term, and the artificial viscosity coefficient is constructed to be proportional to an entropy residual of the system. For scalar conservation laws, the FE residual is one of the entropy residuals. Therefore, the residual-based artificial viscosity method can be obtained by disregarding the GLS terms in the stabilization method of \citep{Szepessy_1989}. In \citep{Nazarov_2013}, we proved that the FE approximation of residual-based artificial viscosity method applied to scalar conservation laws converges to the unique entropy solution, and we extended the method to solve more general systems in the framework of DG, spectral elements, finite differences, and radial basis functions (RBF), see \eg \citep{Nazarov_Hoffman_2013, Marras_2015, Lu_2019, Striernstrom_2021, Tominec_2021}. For other approaches to constructing the nonlinear artificial viscosity methods, we refer to the work of \citep{Basting_2017, Kuzmin_2020, Mabuza_2020}, where the stabilization is constructed based on a smoothness indicator of the solution.

In this work, we propose a new residual-based shock-capturing method for FE approximation of the MHD system. The method is explicit in time and works for high order polynomial degrees. The divergence-free constraint is obtained using several techniques, including the projection of $\bB$ into a divergence-free space  \citep{Brackbill_Barnes_1980} and the hyperbolic cleaning \citep{Dedner_et_al_2002}. We should emphasize that positivity and invariant domain preservation of the method as in \citep{Guermond_et_al_2018}  is left beyond the scope of this paper.

The paper is organized as follows: in Section~\ref{Sec:prelim}, we introduce important notations and terms that are used in the paper, such as FE spaces, meshes, and matrices. In Section~\ref{Sec:Equation}, we introduce the ideal MHD equations, eigenvalues, and the divergence cleaning techniques. In Section~\ref{sec:nv}, we present a viscous regularization of the MHD equations using a novel residual-based shock-capturing method. Finally, in Section~\ref{sec:num}, we solve several well-known benchmark problems for ideal MHD. The method is shown to have an optimal convergence rate for smooth problems and captures discontinuous waves for non-smooth tests. Section~\ref{sec:conclusion} is the conclusion.

\section{Continuous Galerkin finite element method}\label{Sec:prelim}
In this section, we introduce the finite element spaces and notations that are used in this work. 

Let us consider an open bounded domain $\Omega \in \mR^d$, $d=1,2,3$ with boundary $\Gamma$.  We denote $\calT_h$ a finite partition of $\Omega$ into disjoint elements $K$ such that no vertex of any element is placed on the interior of an edge of another element. The union of all the closure of elements $K$ constitutes $\overline\Omega$, where $\overline\Omega$ denotes the closure of $\Omega$.

We seek finite element solution in Lagrange polynomial spaces
\begin{equation}\label{eq:Qh}
  \calQ_h = \{v(x): v \in \calC^0(\Omega), v\vert_K \in \polP_k,\;\forall K \in \calT_h\},
\end{equation}
where $\polP_k$ is a space of polynomials of at most $k$ degrees. Let $\Phi_i$ be the $k$-degree Lagrange basis function in $\calQ_h$ which takes value 1 at node $i$, and 0 at any other nodes. The functions $\{\Phi_i\}_{i=1}^{N_h}$ form a basis for $\calQ_h$, where $N_h$ is the total number of nodes in $\calQ_h$. We define $\calI(i)$ the set of all nodal points contained within the support of $\Phi_i$. Furthermore, we use the following Hilbert space inner products,
\[
(v,w) = \int_{\Omega} vw\ud x, \quad (\GRAD v, \GRAD w) = \int_{\Omega} \sum_{i=1}^d \frac{\p v}{\p x_i}\frac{\p w}{\p x_i} \ud x.
\]

We will use a mesh-function $h_h(\bx)\in \calQ_h$ which we compute using the following projection: find $h_h(\bx)\in \calQ_h$ such that
\begin{equation}\label{eq:h}
  \sum_{K\in \calT_h} \int_K h_h v \ud \bx = \sum_{K\in \calT_h} \int_K \frac{h_K}{k} v \ud \bx,
  \quad \forall v\in \calQ_h,
\end{equation}
where $h_K$ is the circumradius of the element $K$. 

\section{Governing equations of MHD}\label{Sec:Equation}
Let us denote the time interval $[0, \widehat t \,\,]$, where $\widehat t$ is the final time. For all $(\bx,t) \in D:= \Omega\times[0,\widehat t \,\,]$, we seek solution of the MHD equations $\bsfU:=(\rho, \bbm, E, \bB),\;\rho, E: D\to\mR$,\;$\bbm,\bB:D\to\mR^d$, where $\rho$ is the mass density, $\bbm$ is the momentum, $E$ is the total energy, $\bB$ is the magnetic field. The ideal MHD equations in conservative form with appropriate boundary conditions are defined as follows:
\begin{equation}\label{eq:mhd1}
\p_t \bsfU + \DIV \bsfF_{\calE}(\bsfU) + \DIV \bsfF_{\calB}(\bsfU) = 
0,
\end{equation}
where the nonlinear tensor fluxes $\bsfF_{\calE}(\bsfU), \bsfF_{\calB}(\bsfU)$ are defined as
\begin{equation}\label{eq:mhd2}
	\bsfF_{\calE}:=
	\begin{pmatrix}
	\bbm \\
	\bbm \otimes \bu + p\polI \\
	\bu (E+p) \\
	0
	\end{pmatrix},
	\quad
	\bsfF_{\calB}:=
	\begin{pmatrix}
	0 \\
	-\bbetaa\\
	-\bu \SCAL \bbetaa\\
	\bu \otimes \bB - \bB \otimes \bu \\
    \end{pmatrix},
\end{equation}
with $\bu:=\bbm / \rho$ being the velocity field, and $\bbetaa$ being the Maxwell stress tensor:
\[
\bbetaa = \Big( -\frac12(\bB\SCAL\bB)\polI + \bB\otimes\bB \Big),
\]
where $\polI$ denotes the identity matrix of size $d\times d$. The combination of hydrodynamics and electromagnetism nature of the MHD system can be seen from the MHD equations in the form of \eqref{eq:mhd1}, through the separated fluxes $\bsfF_{\calE}(\bsfU)$ and $\bsfF_{\calB}(\bsfU)$. The thermodynamic pressure $p$ is computed from the equation of state 
\begin{equation}\label{eq:equation_of_state}
p = (\gamma-1) e, 
\end{equation} 
where $e$ is the internal energy and $\gamma$ is the adiabatic gas constant. The internal energy is computed as 
\[
e = E - \frac12 \rho (\bu \SCAL \bu) - \frac12 (\bB \SCAL \bB).
\]
In addition, the magnetic field should satisfy the following divergence free constraint:
\begin{equation}\label{eq:div0}
\DIV \bB = 0.
\end{equation}
Note that the above equations are in non-dimensional form, as suggested by \cite[Sec.~3.1]{Sitaraman_Raja_2013}. 

\subsection{Eigenvalues of the inviscid MHD}
We are going to use the maximum wave speeds to construct our stabilization terms. A good approximation of the local wave speeds is obtained by eigenvalues of the inviscid MHD equations. It is well-known that eight eigenvalues are corresponding to the eight elementary waves for the ideal MHD equations, see \eg \cite{Powell_et_al_1991, Barth_1999, Rossmanith_2006}:
let $\be \in \mR^d$ be a direction vector, $a:=\gamma p/\rho$ be the speed of sound, 
$b := \bB\SCAL\be/\rho^\frac12$ and 
\[
c_{f,s}^2 = \frac12
\Big(
a^2 + \frac{\bB \SCAL \bB}{\rho}
\Big) \pm 
\frac12 
\Bigg(
\Big(
a^2 + \frac{\bB \SCAL \bB}{\rho}
\Big)^2 
- 4a^2 b^2
\Bigg)^\frac12.
\]
The eigenvalues from the smallest to the largest are 
\begin{equation}\label{eq:mhd_eigenvalues}
\lambda_{1,8} = \bu\SCAL\be \mp c_{f}, \quad
\lambda_{2,7} = \bu\SCAL\be \mp b, \quad 
\lambda_{3,6} = \bu\SCAL\be \mp c_{s}, \quad
\lambda_{4,5} = \bu\SCAL\be,
\end{equation}
where the subscripts $_{f}$, $_{s}$ and $_{a}$ correspond to fast and slow magnetosonic waves, and the Alf\'en waves. The characteristics $\lambda_{4,5}$ correspond to entropy and divergence waves. It is clear to see that the first and eighth eigenvalues are the ones of interest, which represent the fastest moving waves. Later in Section~\ref{sec:nv}, we will use the maximum wave speed
\[
\lambda_{\max} := \max_{i = 1,\ldots,8} \vert \lambda_i \vert
\]
in the construction of the nonlinear viscosity term.

\subsection{Divergence cleaning techniques}
To satisfy the magnetic divergence-free condition \eqref{eq:div0}, we have verified that our proposed nonlinear stabilization can be used along with popular techniques for divergence cleaning: the elliptic correction/projection method by \cite{Brackbill_Barnes_1980}, and the hyperbolic correction by \cite{Dedner_et_al_2002,Tricco_2016}. A numerical comparison between the mentioned techniques is included in Section \ref{sec:numerical_div_cleaning_comparison}. 

\subsubsection{Elliptic correction/Projection method}\label{sec:projection_cleaning}
In each time step, instead of using the magnetic solution $\bB$ from the system, $\bB'$ the projection of $\bB$ onto a magnetic divergence-free space is used. \cite{Brackbill_Barnes_1980} proposed using
\[
\bB'=\bB-\nabla \Psi_p,
\]
where $\Psi_p$ is the solution to the Poisson equation $\nabla^2\Psi_p - \DIV\bB=0$. This yields $\DIV \bB' = 0$ as the non-physical part $\nabla \Psi_p$ is removed from the numerical solution $\bB$. After each correction in the magnetic field $\bB$, the dependent variables: pressure $p$, temperature $T$, energy $E$ and the entropy functions are updated accordingly to ensure consistency of the discrete solution. Contrary to the seemingly ad-hoc nature of this method, \cite{Toth_2000} shows that the projection method preserves conservative properties and accuracy of the base scheme. 

\textbf{Pseudo time-stepping}. Numerically minimizing $\vert\DIV \bB'\vert$ by the projection method is considered computationally expensive. A more affordable alternative, proposed by \cite{Hayashi_2005}, is solving the following time-dependent equation for $\Psi_p$,
\begin{equation}\label{eq:pseudo_timestepping}
\partial_{\widetilde{t}}\Psi_p-\nabla^2\Psi_p + \DIV\bB = 0,
\end{equation}
where the variable $\widetilde t$ is a separate time dimension. At equilibrium, the solution to \eqref{eq:pseudo_timestepping} is the one to the Poisson equation. A semi-discretization of \eqref{eq:pseudo_timestepping} is given by
\[
\frac{\Psi_p^{\ell+1}-\Psi_p^{\ell}}{\widetilde\tau}-\nabla^2\Psi_p^{\ell+1} + \DIV \bB^{\ell} = 0,\quad \ell=0,1,2,\dots,
\]
where $\widetilde\tau$ is the pseudo step size. After each pseudo time step, the magnetic solution $\bB$ is corrected $\bB^{\ell+1}=\bB^{\ell}-\nabla \Psi_p^{\ell+1}$. The divergence $\DIV \bB'$ is gradually decreased. One can think of this technique as applying the elliptic correction multiple times. The computational cost depends on how many pseudo steps are performed.

\subsubsection{Hyperbolic correction}\label{sec:hyperbolic_cleaning} An auxiliary scalar variable $\Psi_l$ is added to the system to transport and suppress the divergence violation on the magnetic field. We add a fifth equation to the MHD system
\[
\partial_t\Psi_l + c_h \DIV \bB = - \frac{c_rc_h}{h_h}\Psi_l,
\]
where $c_r = 0.3$ in 2D or $c_r = 1.0$ in 3D \citep{Tricco_2016}, and $h_h(\bx)$ denotes the mesh function computed in \eqref{eq:h}. The coefficient $c_h$ is set to be the maximum wave speed $c_h$ =  $\max_{j=1,...,N_h}\left(\vert\lambda_1(\bN_j)\vert, \vert\lambda_8(\bN_j)\vert\right)$. In addition, a source term is added to the right hand side of \eqref{eq:mhd1},
\begin{equation}
\p_t \bsfU + \DIV \bsfF_{\calE}(\bsfU) + \DIV \bsfF_{\calB}(\bsfU) = \bsfG_{\text{GLM}}(\bsfU),
\end{equation}
where
\[
\bsfG_{\text{GLM}}(\bsfU) = 
-\begin{pmatrix}
0 \\
0 \\
\bB\SCAL \nabla \Psi_l \\
\DIV  (\Psi_l \polI) \\
0
\end{pmatrix}.
\]


\section{Nonlinear viscosity method for MHD}
\label{sec:nv}
The aim of this paper is to solve the inviscid MHD system \eqref{eq:mhd1} using continuous finite element approximations. We use the following two vector valued spaces:
\begin{equation}
	\bcalV_h := [\calQ_h]^{d}, \quad \bcalW_h := \calQ_h \CROSS \bcalV_h \CROSS \calQ_h \CROSS \bcalV_h, 
\end{equation}
where $\calQ_h$ is defined in \eqref{eq:Qh}. Then, we formulate the finite element approximation of the MHD system \eqref{eq:mhd1} as follows: 
find $\bsfU_h(t) :=(\rho_h, \bbm_h, E_h, \bB_h)^\top \in \calC^1([0,\widehat t \,\,]; \bcalW_h)$ such that
\begin{equation}\label{eq:ODE_no_viscous}
\begin{aligned}
\big(\p_t \bsfU_h, \bsfV_h \big) 
+ & \big( \DIV \bsfF_{\calE}(\bsfU_h), \bsfV_h \big) 
+  \big( \DIV \bsfF_{\calB}(\bsfU_h), \bsfV_h \big)
= 0,
\quad \forall\,\, \bsfV_h \in \bcalW_h.
\end{aligned}
\end{equation} 
Here the inner products are computed exactly using appropriate numerical quadrature rules. The finite element approximation of the flux terms in \eqref{eq:ODE_no_viscous} is equivalent to central difference schemes, therefore produces spurious oscillations and cannot be used for strong discontinuities such as shocks and contact lines. Below, we introduce a nonlinear artificial viscosity approach to fix this issue.

\subsection{New class of viscous regularization of the MHD equations}
In this section, we propose using the following vanishing viscosity regularization of the MHD system. The new regularization is obtained by adding an elliptic term $\DIV \bsfF_{\calV}(\bsfU)$ into the MHD system, where $\bsfF_{\calV}(\bsfU)$ denotes a viscous flux tensor.

As a heuristic argument to construct $\bsfF_{\calV}(\bsfU)$, we start from the viscous flux of the resistive MHD equations, see e.g., \citep{Dumbser_2009},
\begin{equation}\label{eq:resistive_mhd_flux}
    \bsfF_{\calV}^0(\bsfU_h):=
    \begin{pmatrix}
    0 \\
    \btau_h \\
    \bu_h \SCAL \btau_h + \kappa_h \GRAD T_h + \eta_h\bB_h \SCAL \big( \GRAD \bB_h - \GRAD \bB_h^\top \big)\\
    \eta_h \big( \GRAD \bB_h - \GRAD \bB_h^\top \big)
    \end{pmatrix}.
\end{equation}
where the temperature $T_h$ is a function in $\calQ_h$. The viscous shear stress tensor $\btau_h$ is defined as
\[
\btau_h = \mu_h \Big( \GRAD\bu_h + \GRAD\bu_h^\top\Big) + \lambda_h \DIV\bu_h\polI,
\]
where $\mu_h \ge 0$ is the artificial dynamic viscosity, $\lambda_h$ is the bulk viscosity of undetermined sign such that $\lambda_h + 2\mu_h > 0$, $\kappa_h \ge 0$ is the artificial heat conduction, $\eta_h \ge 0$ is the artificial electric resistivity of the medium. Note that the electric resistivity $\eta_h$ has the same dimension as the dynamic viscosity coefficient $\mu_h$. We point out that the temperature and velocity field are computed at the nodal points: $T_h(\bN_i) = p_h(\bN_i)/\rho_h(\bN_i)$, $\bu_h(\bN_i) = \bbm_h(\bN_i)/\rho_h(\bN_i)$, $i\in \{1,\ldots,N_h\}$.

By considering zero magnetic field, i.e., $\bB \equiv 0$ and thus $\bsfF_{\calB}(\bsfU) \equiv 0$, we obtain from \eqref{eq:mhd1} the compressible Euler equations. For the compressible Euler equations, a frequent choice of viscous regularization is the Navier-Stokes flux, see e.g., \citep{Nazarov_Larcher_2017}, which coincides with $\bsfF_{\calV}^0(\bsfU_h)$ when $\bB \equiv 0$. The Navier-Stokes flux does not add any viscosity term to the mass equation. Therefore, positivity of density can be violated \citep{Guermond_2014, Nazarov_Larcher_2017}. In \cite{Guermond_2014}, the authors have shown the importance of regularizing the mass equation and have proven that it is a key argument for keeping density and internal energy positive, and for the minimum entropy principle. Numerical validation of the viscous regularization compressible Euler equation including the mass equation in the context of high-order finite element approximation is reported in \citep{Nazarov_Larcher_2017}.
They show that the new viscous flux performs better than the traditional Navier-Stokes flux in capturing complex phenomena, e.g., shocks, rarefaction, contact lines, while maintaining high order accuracy.

For MHD, we add a corresponding elliptic term scaling by a positive number $\nu_h$ to regularize the mass conservation. Therefore, the obtained viscous flux is a slightly modified Navier-Stokes viscous fluxes for the flow part and resistivity magnetic fluxes for the magnetic part:
\begin{equation}\label{eq:mhd3}
    \bsfF_{\calV}(\bsfU_h):=
    \begin{pmatrix}
    \nu_h\GRAD\rho_h \\
    \btau_h \\
    \bu_h \SCAL \btau_h + \kappa_h \GRAD T_h + \eta_h \bB_h \SCAL \big( \GRAD \bB_h - \GRAD \bB_h^\top \big)\\
    \eta_h \big( \GRAD \bB_h - \GRAD \bB_h^\top \big)
    \end{pmatrix}.
\end{equation}
Then, the viscous regularized finite element approximation of the MHD equations is: find $\bsfU^{\e}_h(t) \in \calC^1([0,\widehat t \,\,]; \bcalW_h)$ such that
\begin{equation}\label{eq:ODE}
\begin{aligned}
\big(\p_t \bsfU^{\e}_h, \bsfV_h \big) 
+ & \big( \DIV \bsfF_{\calE}(\bsfU^{\e}_h), \bsfV_h \big) 
+  \big( \DIV \bsfF_{\calB}(\bsfU^{\e}_h), \bsfV_h \big) \\
= &
- \big( \bsfF_{\calV}(\bsfU^{\e}_h), \GRAD \bsfV_h) 
+ \big( \bn\SCAL\bsfF_{\calV}(\bsfU^{\e}_h), \bsfV_h)_{\Gamma},
\quad \forall\,\, \bsfV_h \in \bcalW_h,
\end{aligned}
\end{equation}
where $\bn$ is the normal vector pointing outward at every point on the boundary $\Gamma$. In the rest of the paper, we look for a viscous solution, however, to make the notation easier, we drop the superscript $^{\e}$ sign.

We construct the viscosity coefficients $\nu_h, \kappa_h$, and $\eta_h$, all are from the finite element space $\calQ_h$, and make them proportional to a nonlinear functional such that enough viscosity is added to the system in terms of stability and accuracy. 

Let  
$\mu_h^{L,n} := \sum_{i=1}^{N_h}\mu_i^{L,n} \Phi_i$, 
$\mu_h^{H,n} := \sum_{i=1}^{N_h}\mu_i^{H,n} \Phi_i$, 
and 
$\mu_h^n := \sum_{i=1}^{N_h}\mu_i^{n} \Phi_i$
be the low, high and dynamic viscosity coefficients at time $t=t^n$, and $h_h(\bx) := \sum_{i=1}^{N_h}h_i \Phi_i$.

\subsection{First-order viscosity}
At each node $i \in \{1,2,\dots,N_h \}$, we compute the first order viscosity at time $t=t^n$ as
\[
\mu_i^{L,n} = C_{\max}h_i\lambda^n_{\max,i},
\]
where $\lambda^n_{\max}$ is the largest eigenvalue at time $t=t^n$, given by \eqref{eq:mhd_eigenvalues}, and $C_{\max}$ is a positive parameter. A typical range of $C_{\max}$ is reportedly $[0.15, 0.5]$, see \citep{Nazarov_Hoffman_2013}. In one space dimension, it can shown that $C_{\max} = 0.5$ leads to the first order Lax-Friedrichs scheme. We use $C_{\max} = 0.5$ in all the numerical examples in this article for all polynomial degrees. Note, that the mesh function is already scaled with the polynomial degree in \eqref{eq:h}.

This choice of first-order viscosity is consistent with the CFL condition discussed later. In addition, the largest eigenvalue turned out to be sufficient to construct a stable first-order viscosity for the numerical examples considered in this work. We emphasize that for a positivity preserving first-order scheme, similar to the Euler system, as in \cite{Guermond_Popov_2016}, a more accurate estimate of the maximum wave velocity is required. 

Our goal is to build a high-order stabilized method. We use the first-order scheme as a safety factor and then construct a high-order extension of it.

\subsection{High-order residual-based viscosity}
The idea behind a high order scheme is to build an indicator that tracks strong shocks and discontinuities and activates the first-order method from the section above. At the same time, the indicator should decrease the viscosity where the solution is smooth. For this, we use the finite element residual of the PDE as an indicator, since it vanishes in smooth regions, \ie $\sim \calO (h_h^{k+1})$, and has a jump of size $\sim \calO (h_h^{- 1})$ in discontinuities.

We use the residual as an indicator, therefore, we need to get rid of the corresponding solution unit. Our aim is to find a normalization function such that it transforms any input function to the same arithmetic range, and is locally amplified at the regions where the input function changes drastically.

Consider a continuous scalar function $w:\mR^d\to\mR$ and in particular the case $w \in \calQ_h$. We propose dividing $w$ by the following function to normalize $w$
\begin{equation}\label{eq:rv_normalization}
n(w)_i = \left\vert\bar{S}(w) - \alpha\left(\max_{j\in\calI(i)} w_j-\min_{j\in\calI(i)}w_j\right) \right\vert,
\end{equation}
where $\bar{S}$ indicates the global maximum variance $\bar{S}(w) = \| w - \text{mean}(w)\|_{L^{\infty}(\Omega)}$, and the global constant $\alpha$ is calculated as
\[
\alpha = \frac{C_l \bar{S}(w)}{\max_{\Omega}w-\min_{\Omega}w}
\]
which acts as a unit normalizer of local jumps to the global variance. As a technical detail, notice that when $w$ is a constant uniform function, we have ${\max_{\Omega}w-\min_{\Omega}w}=\max_{j\in\calI(i)} w_j-\min_{j\in\calI(i)}w_j=0$ for all $i\in\{1,2,\dots,N_h\}$. In that case, $n(w)$ reduces to $\left\vert(1-C_l)\bar{S}(w)\right\vert$ which although being zero, we do not suffer from division by zero in the calculation of $n(w)$. In an earlier development, the normalization function $n(w)$ in \citep{Striernstrom_2021} can be obtained by setting $\alpha = 1$ in \eqref{eq:rv_normalization}. We observe that the new formulation \eqref{eq:rv_normalization} deals better with mixed-signed solutions since the two factors in $n(w)$ are of the same unit. To demonstrate this argument, we rewrite the expression \eqref{eq:rv_normalization} as
\begin{equation}\label{eq:rv_normalization2}
n(w)_i = \bar{S}(w)\left(1 - C_l \frac{\max_{j\in\calI(i)}w_j-\min_{j\in\calI(i)}w_j}{\max_{\Omega}w-\min_{\Omega}w} \right).
\end{equation}
The parameter $C_l \in (0, 1)$ allows fine tuning of how sensitive this normalization affects the local discontinuities. Observe that the ratio $(\max_{j\in\calI(i)}w_j-\min_{j\in\calI(i)}w_j)/(\max_{\Omega}w-\min_{\Omega}w)$ is in between the interval $[0,1]$, and is amplified at local solution jumps. If the input function $w$ holds both positive and negative parts, and if $\alpha=1$, it is possible that $n(w)_i = 0$ at the regions where the local jumps incidentally matches the global variance. However, once the local jump bypasses the global variance, $n(w)$ can become uncontrollably big, which opposes to our aim that $n(w)$ should be small at such circumstances. On the other hand, the formula \eqref{eq:rv_normalization2} does not suffer this issue because the value within the bracket in \eqref{eq:rv_normalization2} is strictly positive. How small $n(w)$ can become is tuned with the parameter $C_l$. For example, if $C_l = 1$, $n(w)$ can become zero at the largest jump of $w$. Lowering $C_l$ from $1$ towards $0$ would make $n(w)$ at jumps of $w$ larger and larger but not exceeding $\bar{S}(w)$. In smooth regions, $n(w)$ turns into $\bar{S}(w)$, which is the classical normalization function of our residual-based viscosity technique, see \cite{Nazarov_2013}.

Next, we construct a normalized finite element residual of the MHD equations as
\begin{equation}\label{eq:normalized_residual}
\widetilde R_h^n := \max\left( \frac{R(\rho^n_h)}{n(\rho^n_h)}, \frac{R(\bbm^n_{h,1})}{n(\bbm^n_{h,1})}, \dots, \frac{R(\bbm^n_{h,d})}{n(\bbm^n_{h,d})},\frac{R(E^n_h)}{n(E^n_h)},\frac{R(\bB^n_{h,1})}{n(\bB^n_{h,1})},\dots,\frac{R(\bB^n_{h,d})}{n(\bB^n_{h,d})}\right),
\end{equation}
where the normalization terms $n(\rho^n_h),$ $ n(\bbm^n_{h,1}),$ $\dots,$ $ n(\bbm^n_{h,d}),$ $n(E_h),$ $n(\bB^n_{h,1}),\dots,n(\bB^n_{h,d})$, as defined above, normalize the residual units and magnitude. It is also possible that the denominators in \eqref{eq:normalized_residual} become close to zero, in the cases when the initial solution or the analytic solution of the corresponding component is a constant or almost a constant. To avoid division by zero, instead of using $R(q_h)/n(q_h)$, where $q_h$ is a scalar function, we use $R(q^n_h)\frac{n(q^n_h)}{n(q^n_h)^2+\epsilon}$. We choose $\epsilon$ to be machine epsilon, i.e. approximately $2.2e^{-16}$.

For each component $q^n_h\in\{\rho^n_h, \bbm^n_{h,1},\dots,\bbm^n_{h,d},E^n_h,\bB^n_{h,1},\dots,\bB^n_{h,d}\}$, the residual $R(q^n_h)$ can be calculate by standard finite element projection,
\begin{equation}\label{eq:residual}
\int_{\Omega} R(q^n_h) \varphi_i \ud x = \int_{\Omega} \vert \textrm{BDF}(q_h)^n+\DIV(f(q^n_h))\vert\, \varphi_i \ud x,
\end{equation}
where $f(q^n_h)$ corresponds to the nonlinear MHD flux in \eqref{eq:ODE_no_viscous} for each component $q^n_h$, and $\textrm{BDF}(q_h)^n$ is a backward difference approximation of the time derivative $\p_tq_h$ at time $t=t^n$. Note that the second-order BDF scheme is sufficient in getting 4th order accurate scheme.

For  $\polP_1,\polP_3$ Lagrange elements, we observe that without destroying robustness, the residual can be calculated much more efficiently by approximating the above projection solution by lumping the mass matrix,
\[
R(q^n_h) \approx \frac{1}{m_i} \int_{\Omega} \vert \textrm{BDF}(q_h)^n+\DIV(f(q^n_h))\vert\, \varphi_i \ud x,
\]
where $m_i = \int_{\Omega} \varphi_i \ud x > 0$ for $\polP_1,\polP_3$ Lagrange polynomials.

Finally, we are now ready to construct the final artificial viscosity on each node $i$ at time $t=t^n$:
\begin{equation}\label{eq:mu}
\mu_{i}^n = \min\left(\mu_i^{L,n}, \mu_i^{H,n}\right),
\end{equation}
where $\mu_i^H = C_Rh_{h,i}^2 R^n_i$, $R^n_i$ is the nodal value of the residual $\widetilde R_h^n$, $C_R$ is a positive parameter. A typical range of $C_{R}$ is reportedly $[0.1, 1]$, see \citep{Nazarov_Hoffman_2013}. For the numerical simulations in this paper, we use $C_R = 1$. In the rest of this paper, we refer to using \eqref{eq:mu} in the formulation of the viscous regularized equation \eqref{eq:ODE} as the ``RV method''.

\subsection{Adaptive time-stepping}
Once the MHD system is discretized in space using the continuous finite element method we obtain the system of ODEs \eqref{eq:ODE}. Let us denote this system as
\[
	\polM \p_t \bsfU_h(t) = \calF(\bsfU_h(t), \bmu_h(t)),
\]
where $\polM\in\mR^{(2d+2)N_h\times (2d+2)N_h}$ is the mass matrix, $\calF(\bsfU_h(t), \bmu_h(t))$ is the right-hand-side function of the system which depends on the solution $\bsfU_h(t)$ and the viscosity vector is $\bmu_h(t) := (\nu_h(t), \mu_h(t), \kappa_h(t), \eta_h(t))$.

Next, discretize this system in time using explicit $r$-stage Runge-Kutta methods:
\begin{equation}\label{eq:rk}
	\bsfU^{n+1}_h = \bsfU^{n}_h + \tau_n (b_1 \bsfK_1 + \ldots + b_r \bsfK_r),
\end{equation}
where $\tau_n := t^{n+1} - t^n$, $\bsfU^{n}_h := \bsfU_h(t^n)$, $b_i$, $i=1,\ldots,r$ are coefficients obtained from the Butcher tableau, and the stage variables $\bsfK_i$ are computed as follows: for the given solution $\bsfU_h^n$ and the viscosity vector $\bmu_h^n:=\bmu_h(t^n)$ at time level $t^n$, set $\bsfW_0:= \bsfU_h^n$, then let $\bsfW_l$ be the solution at the $l$-th stage of the Runge-Kutta method, then compute $\bsfK_l$ by solving the following system:
\[
	\polM \bsfK_l = \bsfF(\bsfW_l(t), \bmu_h^n),
\]
for all $l=0,\ldots, r$. Note that the viscosity coefficients can also be computed on the fly at every Runge-Kutta stage. However, in this work, the viscosity coefficient is constructed from the previous time level and does not change within the Runge-Kutta stages.

The time-step $\tau_n$ is computed adaptively using a CFL condition
\begin{equation}\label{eq:time_step}
\tau_n = \text{CFL} \frac{\min_{\bx\in\Omega} h_h(\bx)}{\max_{j=1,\ldots,N_h}\vert\lambda^n_{\max,j}\vert}.
\end{equation}

\subsection{Boundary conditions}\label{sec:boundary_conditions}
In the following numerical examples we use the following boundary conditions: Dirichlet, Neumann, periodic, and slip or an impermeability boundary condition. The Neumann boundary condition is implemented as an additional boundary integral in the variational form. The periodic boundary condition is implemented by enforcing the same value of degrees of freedom (DOF) of matching nodes on two opposite boundaries in question.  

The Dirichlet and slip boundary conditions are imposed strongly as a correction step after the Runge-Kutta method is applied. Assume that we want to compute solution at time level $t^{n+1}$, $\bsfU_h^{n+1}$, we impose the Dirichlet boundary condition by setting the values of the solution on the boundary nodes $\bN_j$ by its boundary data. 

The slip or an impermeability boundary condition requires $\bu_h^{n+1} \SCAL \bn = 0$ on the boundary in question. This condition is equivalent to $\bbm_h^{n+1} \SCAL \bn = 0$. Then, the boundary value of momentum at the boundary node $\bN_j$, $\bbm_h^{n+1}(\bN_j)$ is replaced by $\bbm_h^{n+1}(\bN_j) - \bbm_h^{n+1}(\bN_j) \SCAL \bn$.

We refer the reader to \cite{Nazarov_Larcher_2017}, where advantages of the strong implementation of the slip boundary condition versus its weak counterpart are discussed, and we refer to \cite{Guermond_et_al_2018}, where the correction technique is discussed for explicit schemes in more detail.

\subsection{Summary of the algorithm}
To summarize, our stabilized finite element solution of the MHD system \eqref{eq:mhd1} can be obtained by the following steps. In each time step, the solution is advanced as
\begin{enumerate}
\item Calculate the residual using \eqref{eq:residual}, combine the low order and high order viscosity to calculate the amount of artificial viscosity using \eqref{eq:mu}.
\item Solve the main linear system \eqref{eq:ODE} where the time derivative is discretized using a Runge-Kutta method \eqref{eq:rk} of order at least $k+1$.
\item If the projection method or the pseudo time-stepping method is used for cleaning divergence of $\bB_h$, use the techniques described in Section \ref{sec:projection_cleaning}. If the hyperbolic correction method in Section \ref{sec:hyperbolic_cleaning} is used, divergence cleaning is already incorporated into the system in the previous step.
\item Apply boundary conditions strongly as described in Section \ref{sec:boundary_conditions}.
\item Compute an approximation of the maximum local velocity using \eqref{eq:mhd_eigenvalues}, and identify the next time step size using \eqref{eq:time_step}.
\end{enumerate}

\section{Numerical examples}\label{sec:num}
In this section, we demonstrate the efficiency of our proposed stabilization method on several well-known benchmark problems.  For the following tests, the residual viscosity parameters are $C_{\max} = 0.5, C_R = 1.0, C_l = 0.4$. We observe that the parameter $C_l$ is not sensitive as in most circumstances, apart from enhanced smoothness at sharp shocks, any choice of $C_l\in(0,1)$ would give satisfactory outcomes. We have used the classical Runge-Kutta 4 for all the tests with time step \eqref{eq:time_step}, where the CFL number is chosen to be $0.3$. In the viscous flux \eqref{eq:mhd3}, the viscosity coefficients $\nu_h, \kappa_h$, and $\eta_h$ are chosen to be $\mu_h$ given by \eqref{eq:mu}.

\subsection{Accuracy test}
The first two test problems where an analytic solution can be derived are presented to show the high order property of the stabilization in the smooth regions.

\subsubsection{Smooth vortex problem, \cite{Wu_et_al_2018}}
\begin{figure}[h!]
     \centering
     \begin{subfigure}{0.49\textwidth}
         \centering
         \includegraphics[width=\textwidth,viewport=100 480 460 780, clip=true]{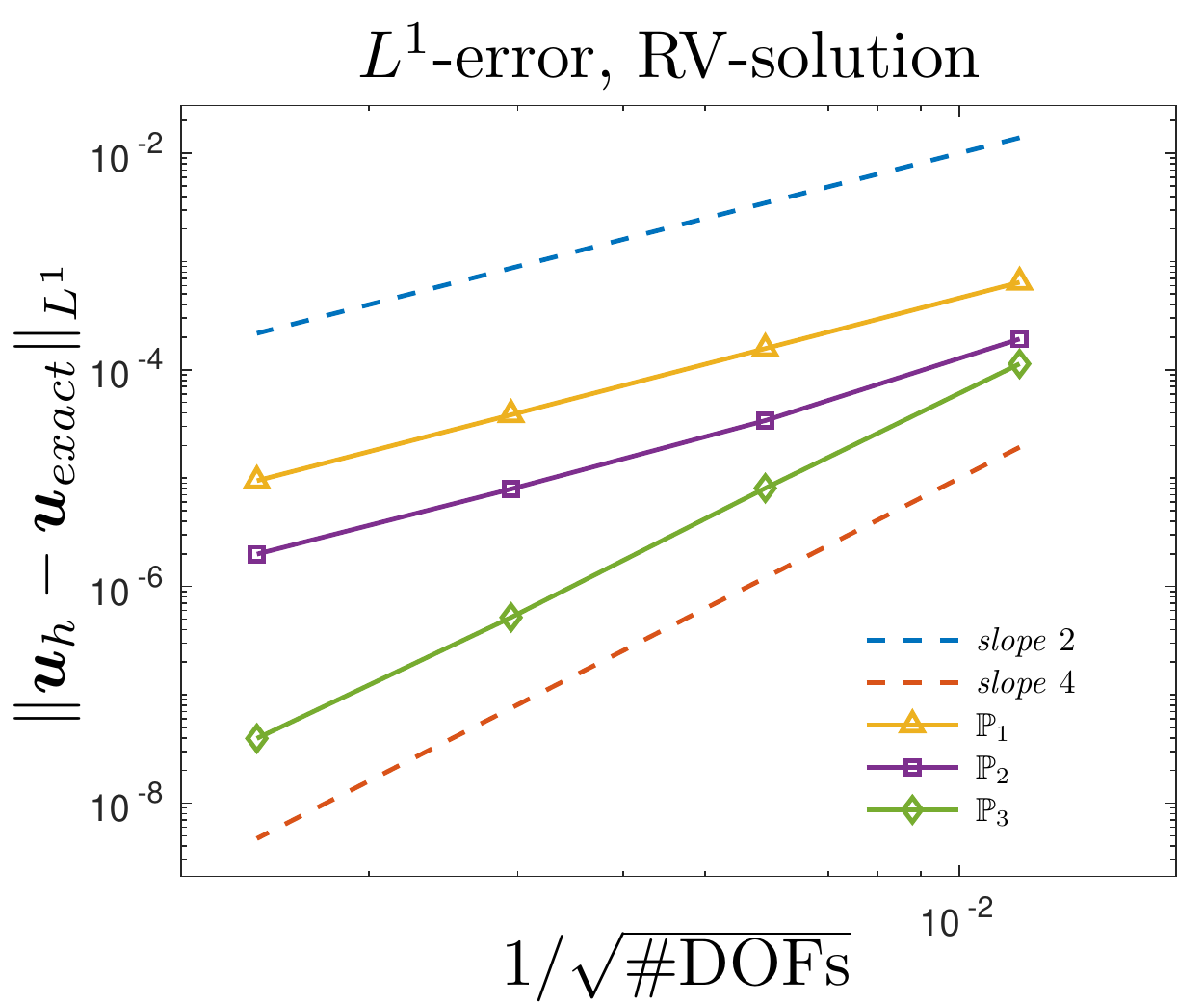}
     \end{subfigure}
     \hfill
     \begin{subfigure}{0.49\textwidth}
         \centering
         \includegraphics[width=\textwidth,viewport=100 480 460 780, clip=true]{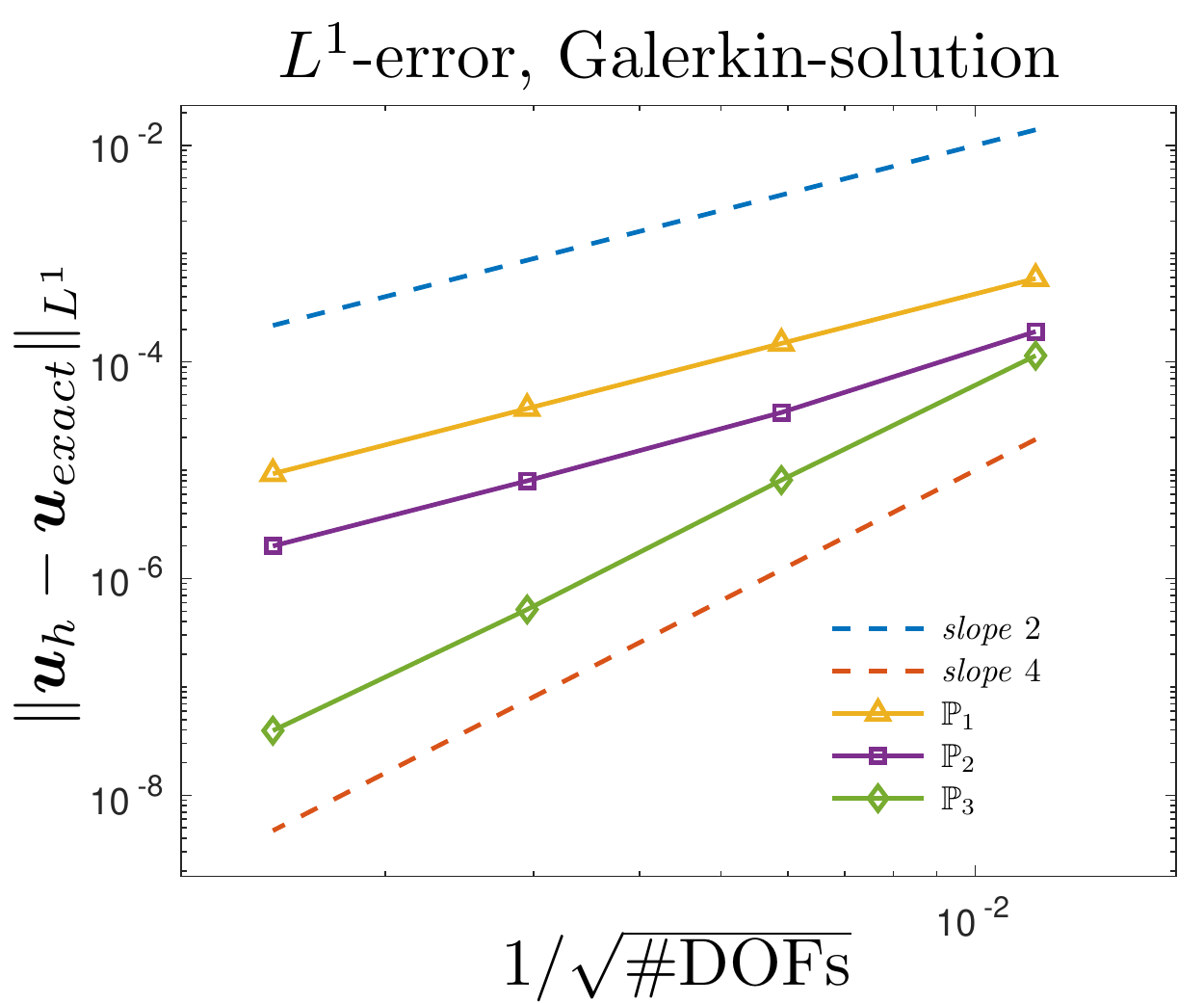}
     \end{subfigure}
     \caption{Accuracy test on the smooth vortex problem: convergence history for velocity $\bu_h$ against the exact solution. The RV method does not destroy the accuracy of high-order solutions.}
     \label{fig:conv_history}
\end{figure}

Consider a periodic smooth vortex problem on a rectangle domain $\Omega = [-10, 10] \times [-10, 10]$. The reference solution is a stationary flow with a vortex perturbation
\[
(\rho(t), \bu(t), p(t), \bB(t)) =(\rho_0, \bu_0+\delta\bu, p_0+\delta p, \bB_0+\delta\bB),
\]
where
\begin{align*}
\rho_0 &= 1, & & \\
\bu_0 &= (1,1), &\delta \bu &= \frac{\mu}{\pi\sqrt 2}e^{(1-r^2)/2}(-r_2,r_1),\\
p_0 &=0, &\delta p &= -\frac{\mu^2(1+r^2)e^{1-r^2}}{8\pi^2}, \\
\bB_0 &=(0.1,0.1), & \delta \bB & = \frac{\mu e^{(1-r^2)/2}}{2\pi}(-r_2,r_1),
\end{align*}
the vortex radius is $r=\sqrt{r_1^2+r_2^2}$, $(r_1,r_2)=(x,y)-\bu_0t$, and the vortex strength is $\mu=5.389489439$. The adiabatic constant is  $\gamma=\frac{5}{3}$. The errors measured at final time $\widehat t=0.05$ using $\polP_1, \polP_2, \polP_3$ elements  are shown in Table \ref{table:convergence_vortex}. In Table \ref{table:convergence_vortex}, we show behavior of both the Galerkin solutions and the RV solutions. For the $\polP_2$ solution, the obtained error is lower than that of the $\polP_1$ solution under the same number of degrees of freedom. However, the convergence rate of the $\polP_2$ is second order, which is suboptimal with regards to theoretical derivations. This is a known issue with continuous Galerkin approximations, which was earlier reported in, e.g., \citep{Nazarov_Larcher_2017}. Figure \ref{fig:conv_history} illustrates and compares the convergence history between the Galerkin solutions and the RV solutions. All the obtained convergence rates show that the RV method does not destroy high-order accuracy of the Galerkin smooth solutions.

\begin{table}
    \centering
    \caption{Smooth vortex problem. L$^1$, L$^2$ error and convergence rates of velocity $\bu$ and magnetic field $\bB$. $\polP_1,\polP_2,\polP_3$ solutions. $\widehat t=0.05$.}
    \label{table:convergence_vortex}
    \begin{tabular*}{\textwidth}{@{\extracolsep{\fill}}ccccccccc@{\extracolsep{\fill}}}
    \toprule%
    \multicolumn{9}{@{}c@{}}{$\polP_1$ solution} \\ \midrule
    \multirow{2}{*}{\#DOFs} &  \multicolumn{4}{@{}c@{}}{RV-solution $\bu$} &  \multicolumn{4}{@{}c@{}}{Galerkin-solution $\bu$} \\\cmidrule{2-5}\cmidrule{6-9}
     {}   &       L$^1$-error   &    Rate   &       L$^2$-error   &    Rate &       L$^1$-error   &    Rate   &       L$^2$-error   &    Rate   \\ \midrule
     7442 & 6.42E-04 & --  &   3.56E-03 & -- &  5.89E-04 &  -- &   3.24E-03 & --  \\ 
    29282 & 1.57E-04 &   2.03 &   8.72E-04 &   2.03  & 1.48E-04 &   1.99 &   8.18E-04 &   1.99  \\ 
   116162 & 3.85E-05 &   2.03 &   2.13E-04 &   2.03 &  3.71E-05 &   2.00 &   2.05E-04 &   2.00 \\ 
   462722 &  9.47E-06 &   2.02 &   5.24E-05 &   2.02 & 9.28E-06 &   2.00 &   5.13E-05 &   2.00 \\ \midrule
    \multirow{2}{*}{\#DOFs} &  \multicolumn{4}{@{}c@{}}{RV-solution $\bB$} &  \multicolumn{4}{@{}c@{}}{Galerkin-solution $\bB$} \\\cmidrule{2-5}\cmidrule{6-9}
     {}   &       L$^1$-error   &    Rate   &       L$^2$-error   &    Rate &       L$^1$-error   &    Rate   &       L$^2$-error   &    Rate   \\ \midrule
     7442 & 2.65E-02 &   --   &   3.02E-02 &  --    &  2.35E-02 &   --   &   2.60E-02 &   -- \\ 
    29282 & 6.49E-03 &   2.06 &   7.36E-03 &   2.06 &  5.90E-03 &   2.01 &   6.56E-03 &   2.01  \\ 
   116162 & 1.57E-03 &   2.06 &   1.76E-03 &   2.07 &  1.48E-03 &   2.01 &   1.64E-03 &   2.01 \\ 
   462722 & 3.82E-04 &   2.04 &   4.27E-04 &   2.05 &  3.70E-04 &   2.01 &   4.11E-04 &   2.01 \\ \mmidrule
   \multicolumn{9}{@{}c@{}}{$\polP_2$ solution} \\ \midrule
    \multirow{2}{*}{\#DOFs} &  \multicolumn{4}{@{}c@{}}{RV-solution $\bu$} &  \multicolumn{4}{@{}c@{}}{Galerkin-solution $\bu$} \\\cmidrule{2-5}\cmidrule{6-9}
     {}   &       L$^1$-error   &    Rate   &       L$^2$-error   &    Rate &       L$^1$-error   &    Rate   &       L$^2$-error   &    Rate   \\ \midrule
     7442 & 1.93E-04 &   -- &   1.14E-03 &   --  &  1.91E-04 &   -- &   1.13E-03 &   -- \\ 
    29282 & 3.40E-05 &   2.50 &   2.01E-04 &   2.50  &  3.41E-05 &   2.48 &   2.02E-04 &   2.49  \\ 
   116162 & 7.94E-06 &   2.10 &   4.67E-05 &   2.10  &  7.98E-06 &   2.10 &   4.70E-05 &   2.10 \\ 
   462722 & 1.99E-06 &   2.00 &   1.17E-05 &   2.00  &  2.00E-06 &   2.00 &   1.17E-05 &   2.00  \\ \midrule
   \multirow{2}{*}{\#DOFs} &  \multicolumn{4}{@{}c@{}}{RV-solution $\bu$} &  \multicolumn{4}{@{}c@{}}{Galerkin-solution $\bu$} \\\cmidrule{2-5}\cmidrule{6-9}
     {}   &       L$^1$-error   &    Rate   &       L$^2$-error   &    Rate &       L$^1$-error   &    Rate   &       L$^2$-error   &    Rate   \\ \midrule
     7442 & 7.96E-03 &   -- &   9.26E-03 &   -- &  7.83E-03 &   -- &   9.08E-03 &   -- \\ 
    29282 & 1.41E-03 &   2.53 &   1.55E-03 &   2.61 &  1.40E-03 &   2.51 &   1.54E-03 &   2.59  \\ 
   116162 & 3.18E-04 &   2.16 &   3.32E-04 &   2.23 &  3.18E-04 &   2.15 &   3.32E-04 &   2.22 \\ 
   462722 & 7.74E-05 &   2.04 &   8.01E-05 &   2.06 &  7.76E-05 &   2.04 &   8.03E-05 &   2.06 \\ \mmidrule
   \multicolumn{9}{@{}c@{}}{$\polP_3$ solution} \\ \midrule
    \multirow{2}{*}{\#DOFs} &  \multicolumn{4}{@{}c@{}}{RV-solution $\bu$} &  \multicolumn{4}{@{}c@{}}{Galerkin-solution $\bu$} \\\cmidrule{2-5}\cmidrule{6-9}
     {}   &       L$^1$-error   &    Rate   &       L$^2$-error   &    Rate &       L$^1$-error   &    Rate   &       L$^2$-error   &    Rate   \\ \midrule
     7442 & 1.14E-04 &   -- &   5.78E-04 &   --  &  1.14E-04 &   -- &   5.87E-04 &   -- \\ 
    29282 & 8.08E-06 &   3.81 &   4.91E-05 &   3.56  &  8.11E-06 &   3.82 &   4.94E-05 &   3.57  \\ 
   116162 & 5.18E-07 &   3.96 &   3.90E-06 &   3.65  &  5.18E-07 &   3.97 &   3.91E-06 &   3.66 \\ 
   462722 & 3.96E-08 &   3.71 &   6.62E-07 &   2.56  &  3.95E-08 &   3.72 &   6.62E-07 &   2.56 \\ \midrule
    \multirow{2}{*}{\#DOFs} &  \multicolumn{4}{@{}c@{}}{RV-solution $\bu$} &  \multicolumn{4}{@{}c@{}}{Galerkin-solution $\bu$} \\\cmidrule{2-5}\cmidrule{6-9}
     {}   &       L$^1$-error   &    Rate   &       L$^2$-error   &    Rate &       L$^1$-error   &    Rate   &       L$^2$-error   &    Rate   \\ \midrule
     7442 & 4.43E-03 &   -- &   4.43E-03 &   -- &  4.44E-03 &   -- &   4.51E-03 &   --  \\ 
    29282 & 2.79E-04 &   4.04 &   2.91E-04 &   3.98 &  2.81E-04 &   4.03 &   2.96E-04 &   3.98  \\ 
   116162 & 1.79E-05 &   3.99 &   2.25E-05 &   3.71 &  1.79E-05 &   3.99 &   2.26E-05 &   3.73  \\ 
   462722 & 1.42E-06 &   3.67 &   2.37E-06 &   3.26 &  1.41E-06 &   3.69 &   2.32E-06 &   3.29  \\
   \botrule
   \end{tabular*}
\end{table}

\subsubsection{Smooth wave propagation, \cite{Wu_et_al_2018}}
In this test case, we consider a rectangle domain $\Omega = [0, 2\pi] \times [0, 2\pi]$. A periodic wave-like solution reads $(\rho(t), \bu(t), p(t), \bB(t))$ $=$ $(1+0.99\sin(x+y-2t),(1,1),1,(0.1,0.1))$. The adiabatic constant is $\gamma=1.4$. The errors measured at final time $\widehat t=0.1$ using $\polP_1, \polP_2, \polP_3$ elements  are shown in Table \ref{table:convergence_constantB}. Again, the obtained convergence rates are optimal with regards to the corresponding Galerkin approximations.

\begin{table}[h!]
    \centering
    \caption{Smooth wave propagation. Errors in density, $\polP_1, \polP_2, \polP_3$ solutions, final time  $\widehat t=0.1$.}
    \label{table:convergence_constantB}
    \begin{tabular*}{\textwidth}{@{\extracolsep{\fill}}ccccccccc@{\extracolsep{\fill}}}
        \toprule%
{}  & &     \#DOFs &       L$^1$-error   &    Rate   &       L$^2$-error   &   Rate   &    L$^\infty$-error   &   Rate  \\ \midrule
\multirow{8}{*}{\rotatebox[origin=c]{90}{$\polP_1$ solution}} & \multirow{4}{*}{\rotatebox[origin=c]{90}{Galerkin}}  &         961 &   6.90E-03 &   --   &   6.87E-03 &   --   &   9.66E-03 &   --   \\
 & &    3721 &   1.73E-03 &   2.05 &   1.72E-03 &   2.04 &   2.42E-03 &   2.04 \\
 & &   14641 &   4.32E-04 &   2.02 &   4.30E-04 &   2.02 &   6.06E-04 &   2.02 \\
 & &   58081 &   1.08E-04 &   2.01 &   1.08E-04 &   2.01 &   1.52E-04 &   2.01  \\  \cmidrule{2-9}
 & \multirow{4}{*}{\rotatebox[origin=c]{90}{RV}}  &        961 &   7.69E-03 &   --   &   7.79E-03 &   --   &   1.23E-02 &   -- \\
 & &    3721 &   1.80E-03 &   2.15 &   1.80E-03 &   2.17 &   2.55E-03 &   2.32 \\
 & &   14641 &   4.41E-04 &   2.05 &   4.40E-04 &   2.05 &   6.17E-04 &   2.07 \\
 & &   58081 &   1.09E-04 &   2.03 &   1.09E-04 &   2.03 &   1.53E-04 &   2.02  \\ \botrule
\multirow{8}{*}{\rotatebox[origin=c]{90}{$\polP_2$ solution}} & \multirow{4}{*}{\rotatebox[origin=c]{90}{Galerkin}}  &       961 &   1.06E-03 &   --   &   1.24E-03 &   --   &   2.54E-03 &   --   \\
 & &    3721 &   1.90E-04 &   2.54 &   2.30E-04 &   2.49 &   5.19E-04 &   2.35 \\
 & &   14641 &   4.29E-05 &   2.17 &   5.14E-05 &   2.19 &   1.34E-04 &   1.98 \\
 & &   58081 &   1.05E-05 &   2.05 &   1.24E-05 &   2.06 &   3.37E-05 &   2.00  \\ \cmidrule{2-9}
 & \multirow{4}{*}{\rotatebox[origin=c]{90}{RV}}  &        961 &   1.44E-03 &   --   &   1.65E-03 &   --   &   2.98E-03 &   --   \\
 & &    3721 &   2.65E-04 &   2.50 &   3.12E-04 &   2.46 &   6.99E-04 &   2.14 \\
 & &   14641 &   4.42E-05 &   2.61 &   5.19E-05 &   2.62 &   1.33E-04 &   2.43 \\
 & &   58081 &   1.05E-05 &   2.08 &   1.24E-05 &   2.07 &   3.35E-05 &   2.00  \\ \botrule
\multirow{8}{*}{\rotatebox[origin=c]{90}{$\polP_3$ solution}} & \multirow{4}{*}{\rotatebox[origin=c]{90}{Galerkin}}  &        961 &   2.33E-04 &   --   &   2.33E-04 &   --   &   7.36E-04 &   --   \\
 & &    3721 &   1.36E-05 &   4.19 &   1.71E-05 &   3.86 &   6.22E-05 &   3.65 \\
 & &   14641 &   7.46E-07 &   4.24 &   1.53E-06 &   3.52 &   5.42E-06 &   3.56 \\
 & &   58081 &   4.61E-08 &   4.04 &   1.67E-07 &   3.21 &   4.69E-07 &   3.55  \\ \cmidrule{2-9}
 & \multirow{4}{*}{\rotatebox[origin=c]{90}{RV}}  &       961 &   2.32E-04 &   --   &   2.33E-04 &   --   &   7.33E-04 &   --   \\
 & &    3721 &   1.36E-05 &   4.19 &   1.70E-05 &   3.86 &   6.21E-05 &   3.65 \\
 & &   14641 &   7.46E-07 &   4.24 &   1.52E-06 &   3.52 &   5.42E-06 &   3.56 \\
 & &   58081 &   4.65E-08 &   4.03 &   1.67E-07 &   3.21 &   4.73E-07 &   3.54  \\
    \botrule
    \end{tabular*}
\end{table}

\subsection{Brio-Wu MHD shock tube problem, \cite{Brio_Wu_1988}}
The RV method shows high-order accuracy for smooth problems considered in above. Now, we start solving problems with strong shocks and discontinuities. 

We want to solve the one-dimensional Brio-Wu MHD shock tube problem, which was first introduced by \cite{Brio_Wu_1988}. This is the one-dimensional Riemann problem in the domain $\Omega = [0, 1]$. The problem is challenging because the solution contains waves of different types, typical for an ideal MHD.

The problem setting is the following: 
\[
(\rho, u, p, B_x, B_y) = 
\begin{cases}
(1,0,1,0.75, 1), &x \in [0,0.5), \\
(0.125,0,0.1,0.75,-1), & x \in [0.5,1.0].
\end{cases}
\]

\begin{figure}[h!]
     \centering
     \begin{subfigure}{0.49\textwidth}
         \centering
         \includegraphics[width=\textwidth,viewport=100 450 480 780, clip=true]{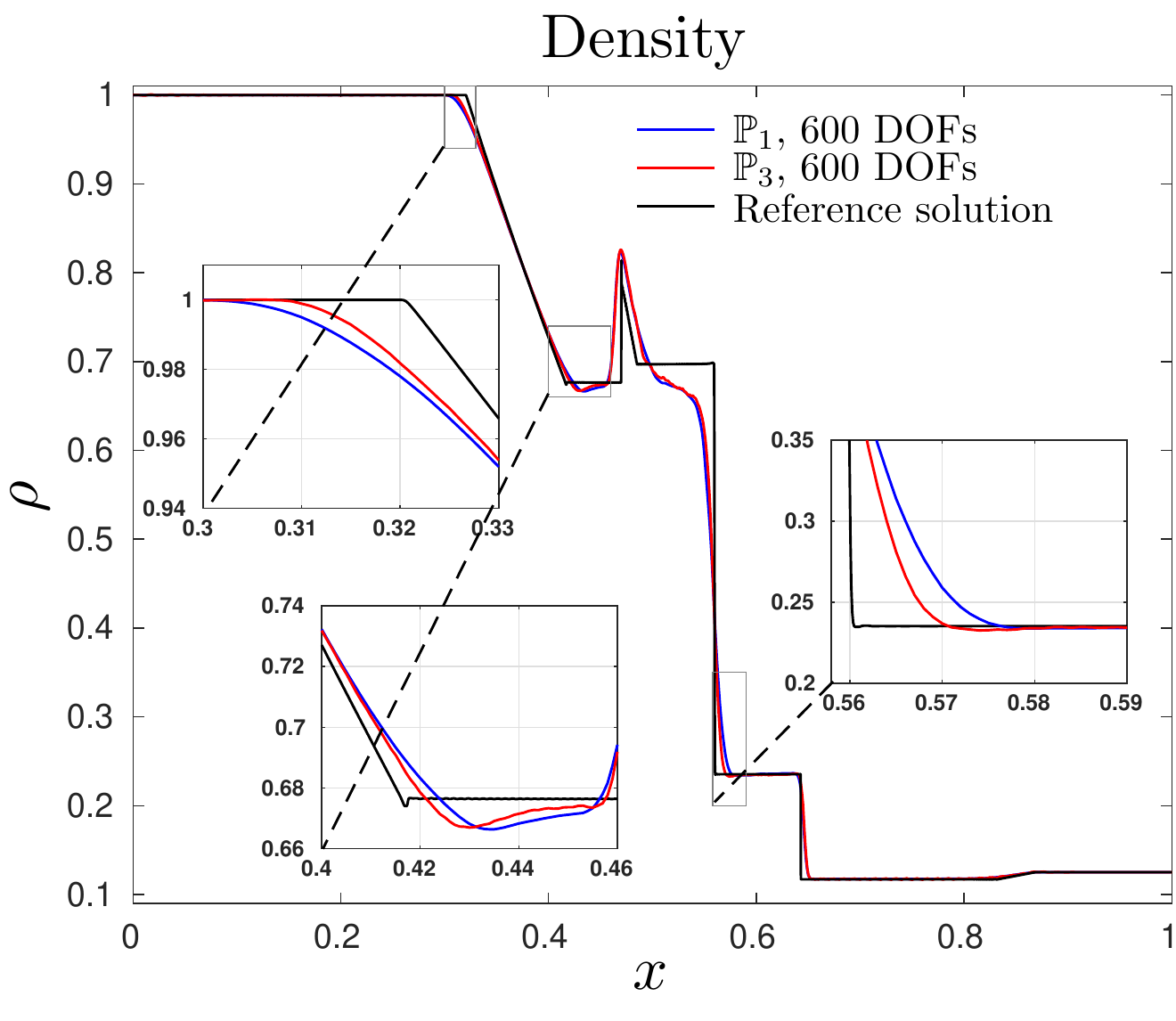}
         \caption{$\polP_1, \polP_3$ solutions with 600 DOFs}
     \end{subfigure}
     \hfill
     \begin{subfigure}{0.49\textwidth}
         \centering
         \includegraphics[width=\textwidth,viewport=100 450 480 780, clip=true]{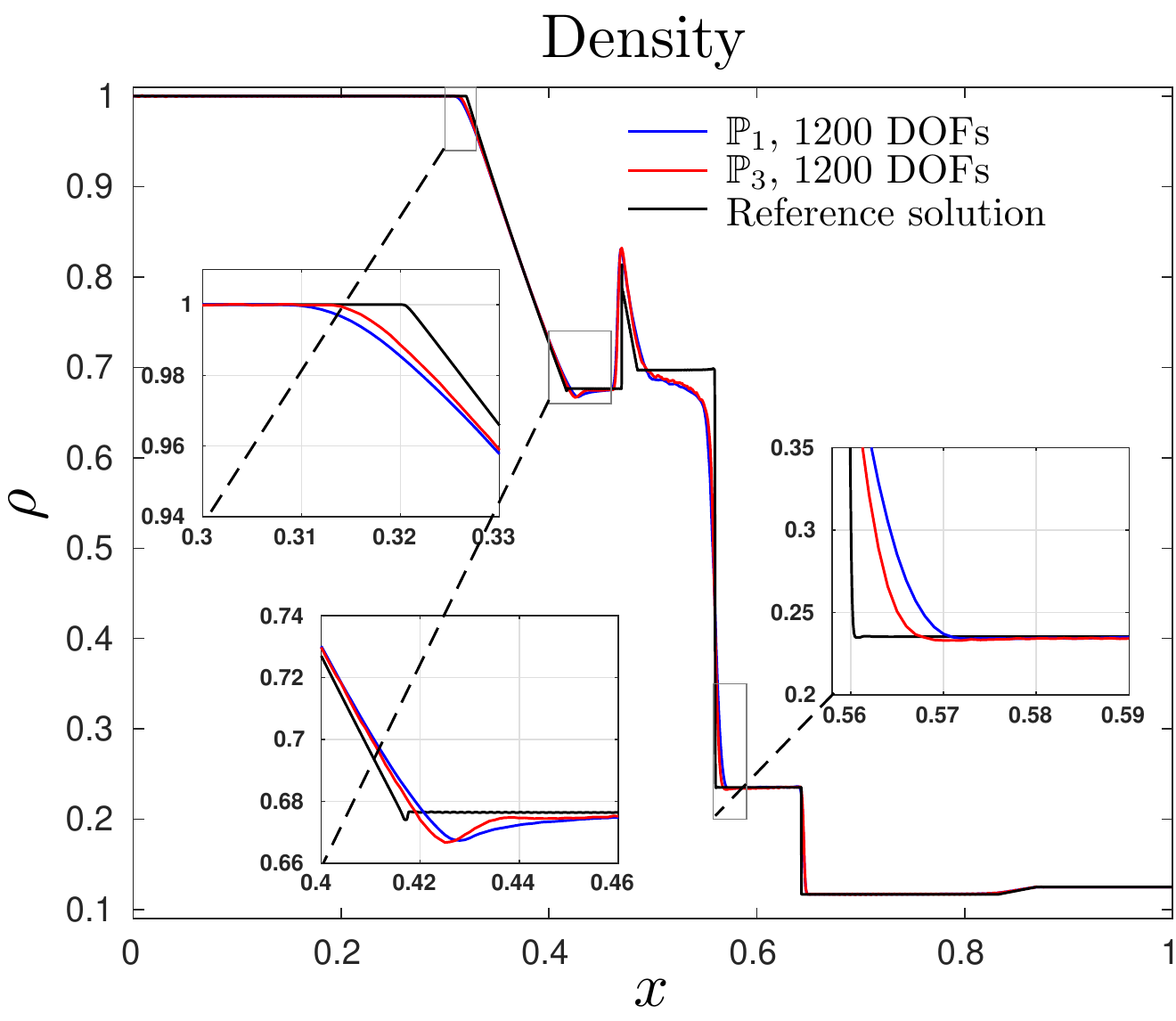}
         \caption{$\polP_1, \polP_3$ solutions with 1200 DOFs}
     \end{subfigure}
     \caption{Comparison of $\polP_1$ and $\polP_3$ RV-solutions to the Brio-Wu problem. The reference solution is obtained from the Athena code \citep{Stone_2008} with 10000 grid points. Under the same number of degrees of freedom, the $\polP_3$ solution captures the compound structure more accurately than the $\polP_1$ solution.}
     \label{fig:brio_wu_P3}
\end{figure}

The adiabatic constant is $\gamma=2$. The importance of this test is to examine whether a MHD solver can represent the shocks, rarefactions, contact lines, and the whole compound structure accurately. The solution at the final time $\hat t = 0.2$ is shown in Figure \ref{fig:brio_wu_P3}. A well-known reference solution is provided by the Athena code \citep{Stone_2008} with 10000 grid points.

The result in Figure \ref{fig:brio_wu_P3} shows that our proposed method can capture the wave structures efficiently and accurately. We also see that $\polP_3$ solution is sharper than $\polP_1$ for the same number of degrees of freedom.  

\subsection{2D Orszag-Tang problem, \citep{Orszag_Tang_1979}}\label{sec:numerical_OT_2D}
\begin{figure}[h!]
     \centering
     \begin{subfigure}{0.49\textwidth}
         \centering
         \includegraphics[width=\textwidth]{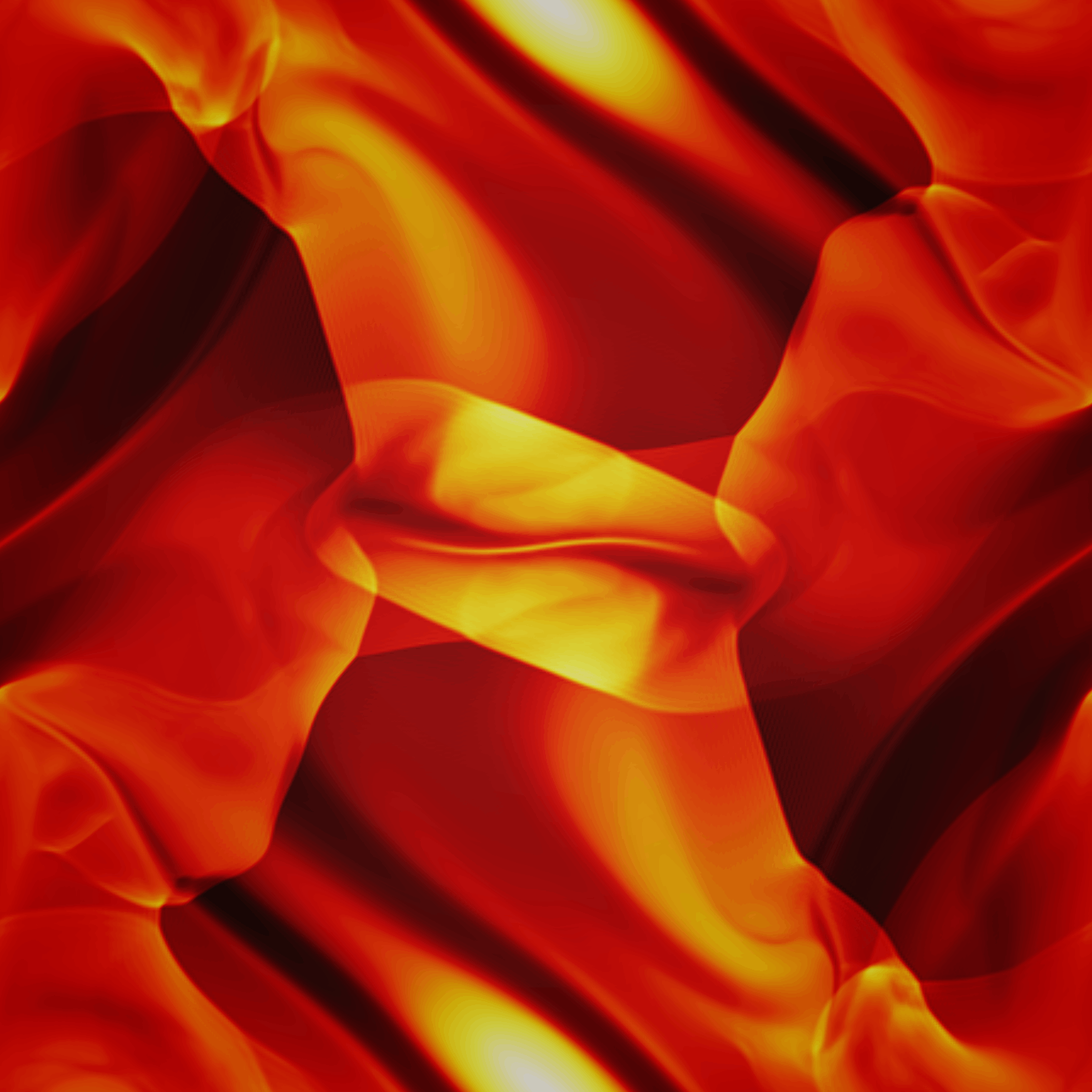}
         \caption{Density at $t=0.5$,\\$\rho_h\in[8.14\text{E-}2,4.95\text{E-}1]$}
     \end{subfigure}
     \hfill
     \begin{subfigure}{0.49\textwidth}
         \centering
         \includegraphics[width=\textwidth]{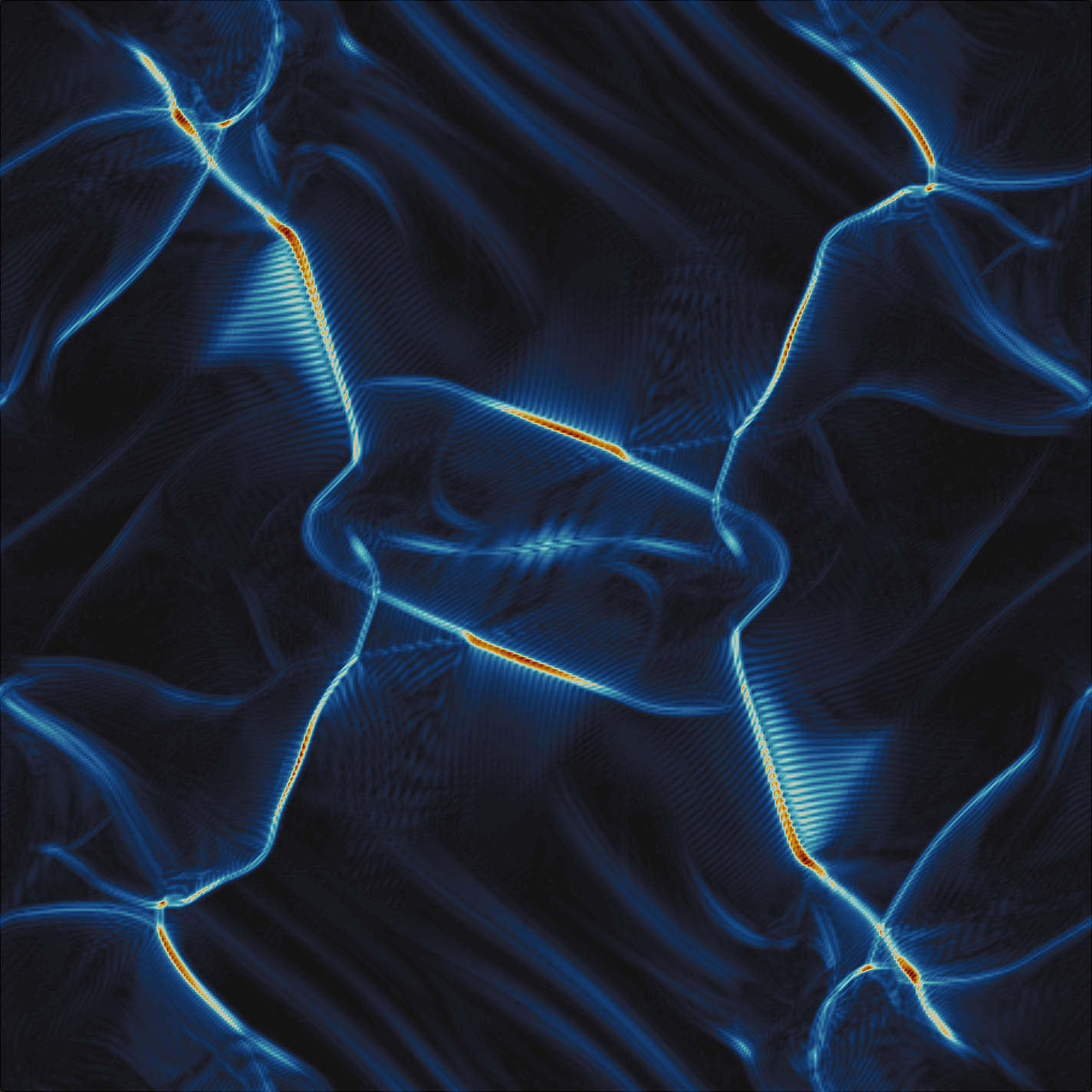}
         \caption{Residual viscosity at $t=0.5$,\\$\mu_h\in[4.09\text{E-}7,4.20\text{E-}4]$}
     \end{subfigure}
     \hfill
     \begin{subfigure}{0.49\textwidth}
         \centering
         \includegraphics[width=\textwidth]{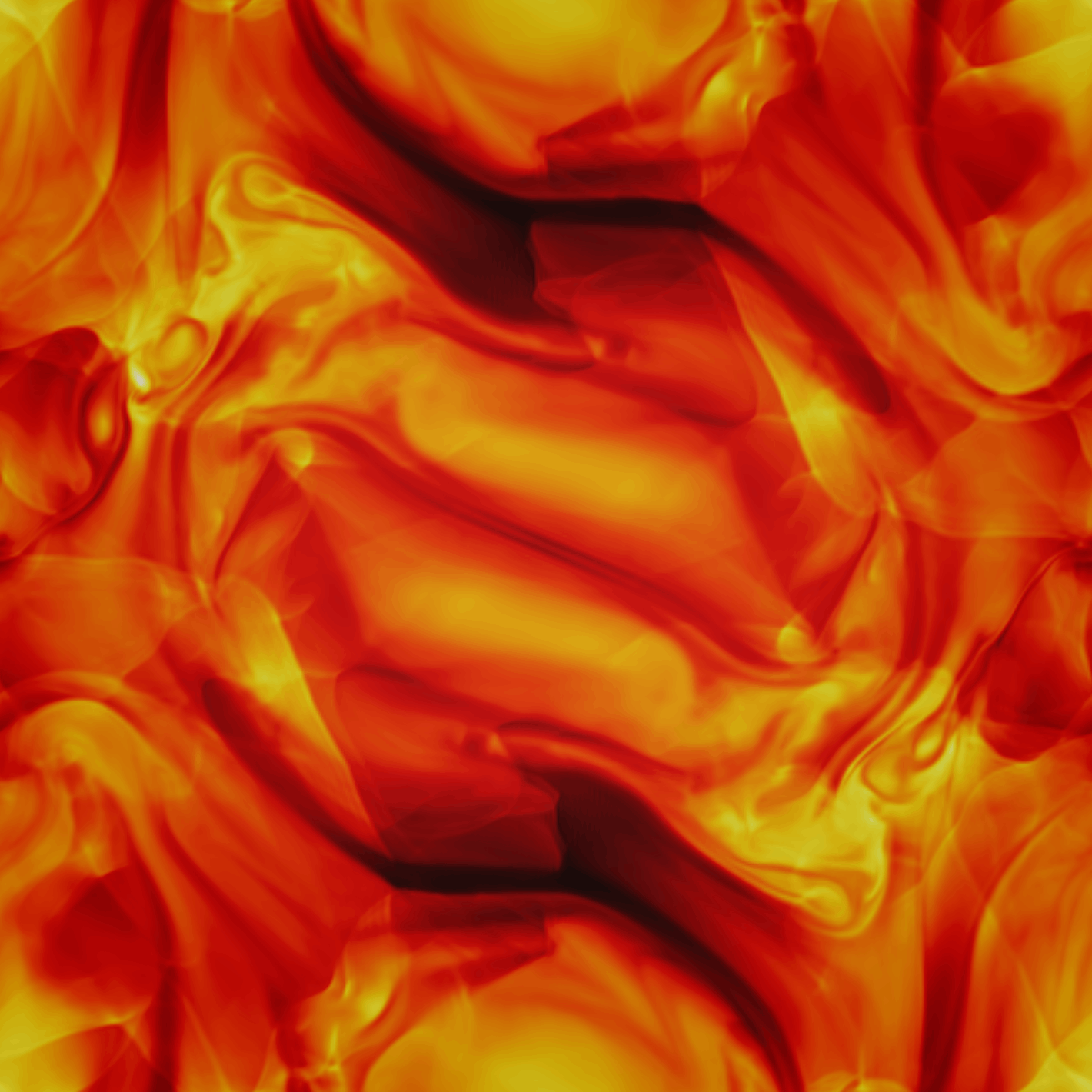}
         \caption{Density at $t=0.5$,\\$\rho_h\in[4.17\text{E-}2,4.41\text{E-}1]$}
     \end{subfigure}
     \hfill
     \begin{subfigure}{0.49\textwidth}
         \centering
         \includegraphics[width=\textwidth]{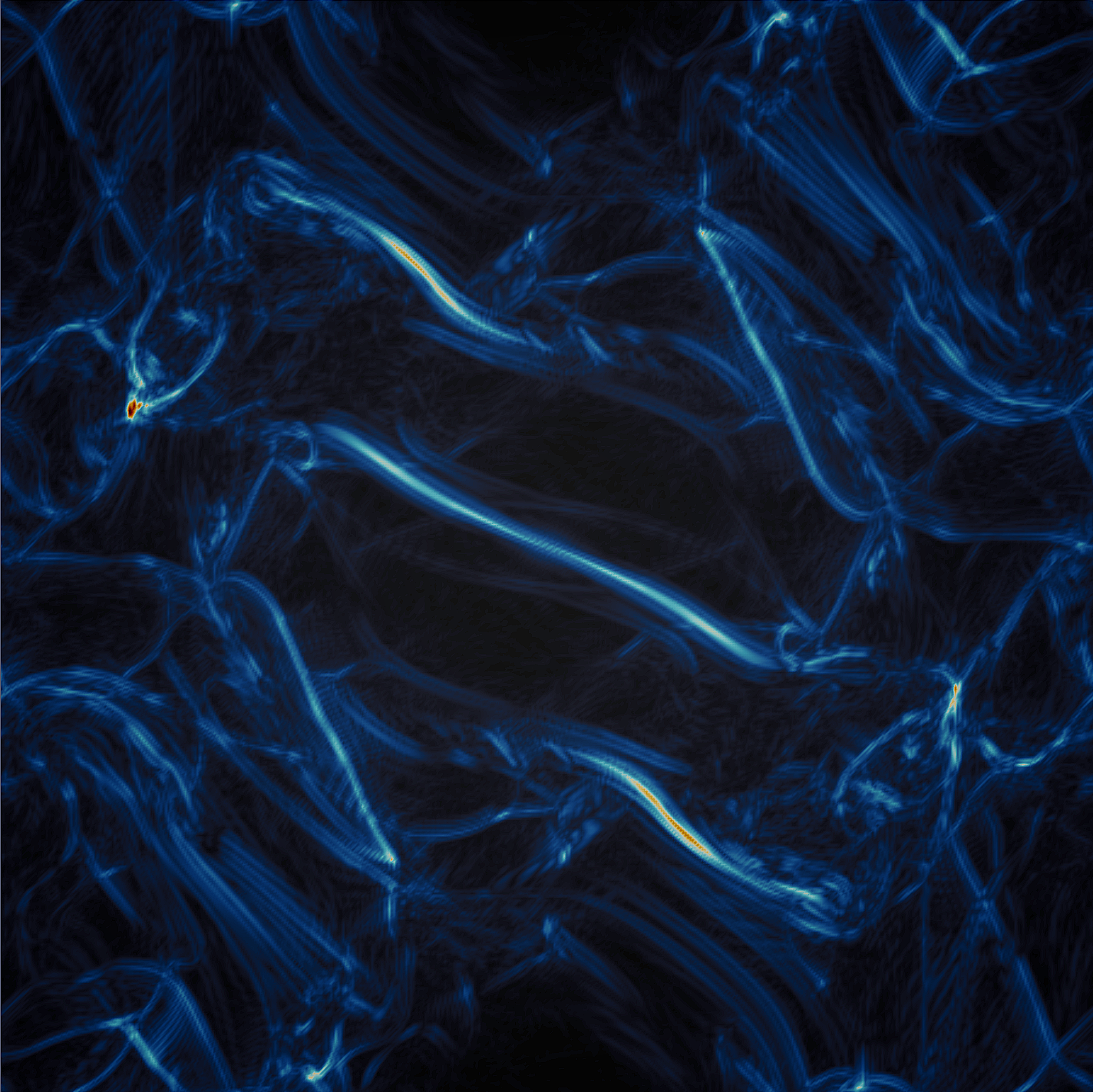}
         \caption{Residual viscosity at $t=1.0$,\\$\mu_h\in[5.54\text{E-}7,4.84\text{E-}4]$}
     \end{subfigure}
     \caption{$\polP_1$ solution of the 2D Orszag-Tang problem, $450\times 450$ nodes. Shocks and turbulence are captured using the residual-based viscosity method.}
     \label{fig:OT_2D_P1}
\end{figure}

\begin{figure}[h!]
     \centering
     \begin{subfigure}{0.49\textwidth}
         \centering
         \includegraphics[width=\textwidth]{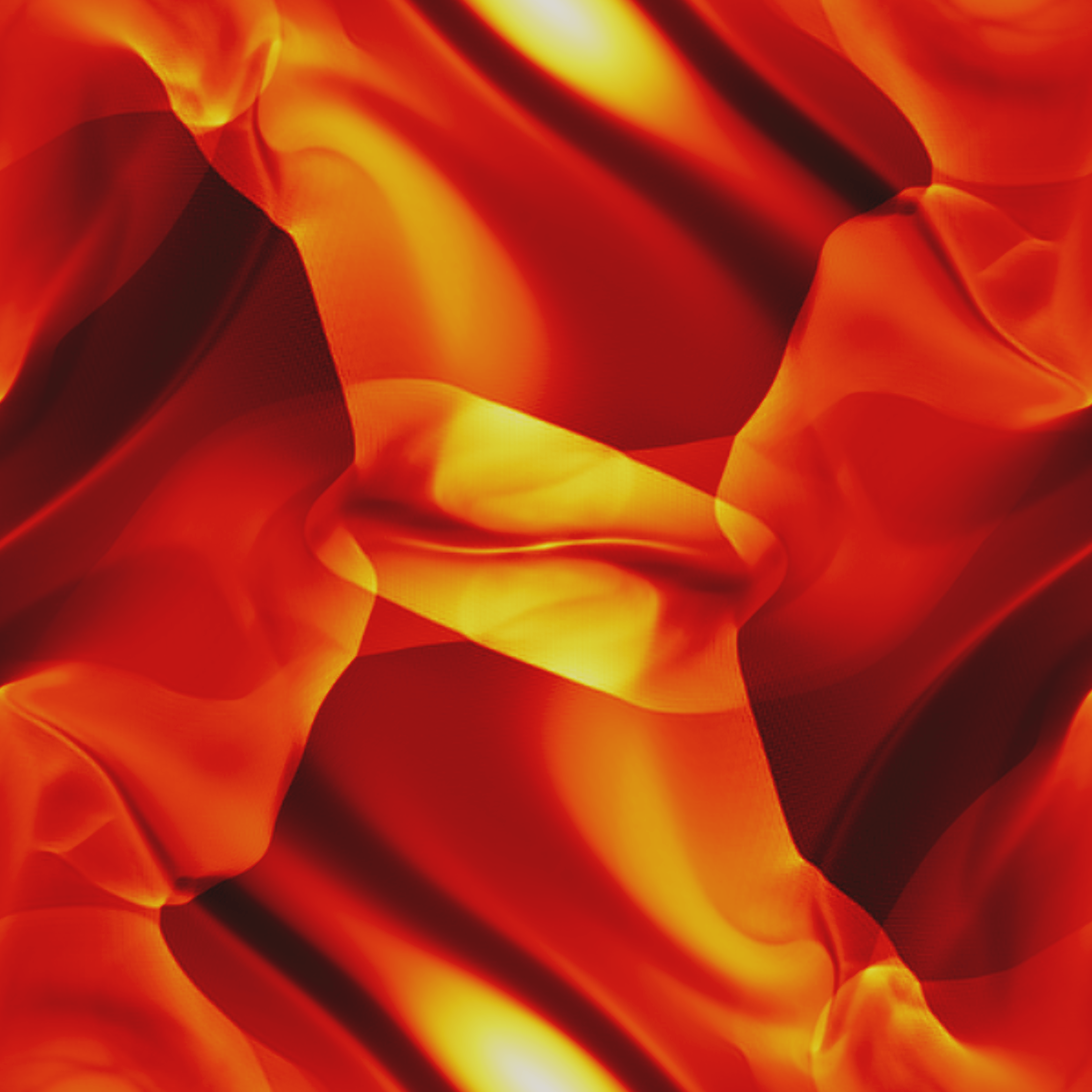}
         \caption{Density at $t=0.5$,\\$\rho_h\in[7.67\text{E-}2,4.96\text{E-}1]$}
     \end{subfigure}
     \hfill
     \begin{subfigure}{0.49\textwidth}
         \centering
         \includegraphics[width=\textwidth]{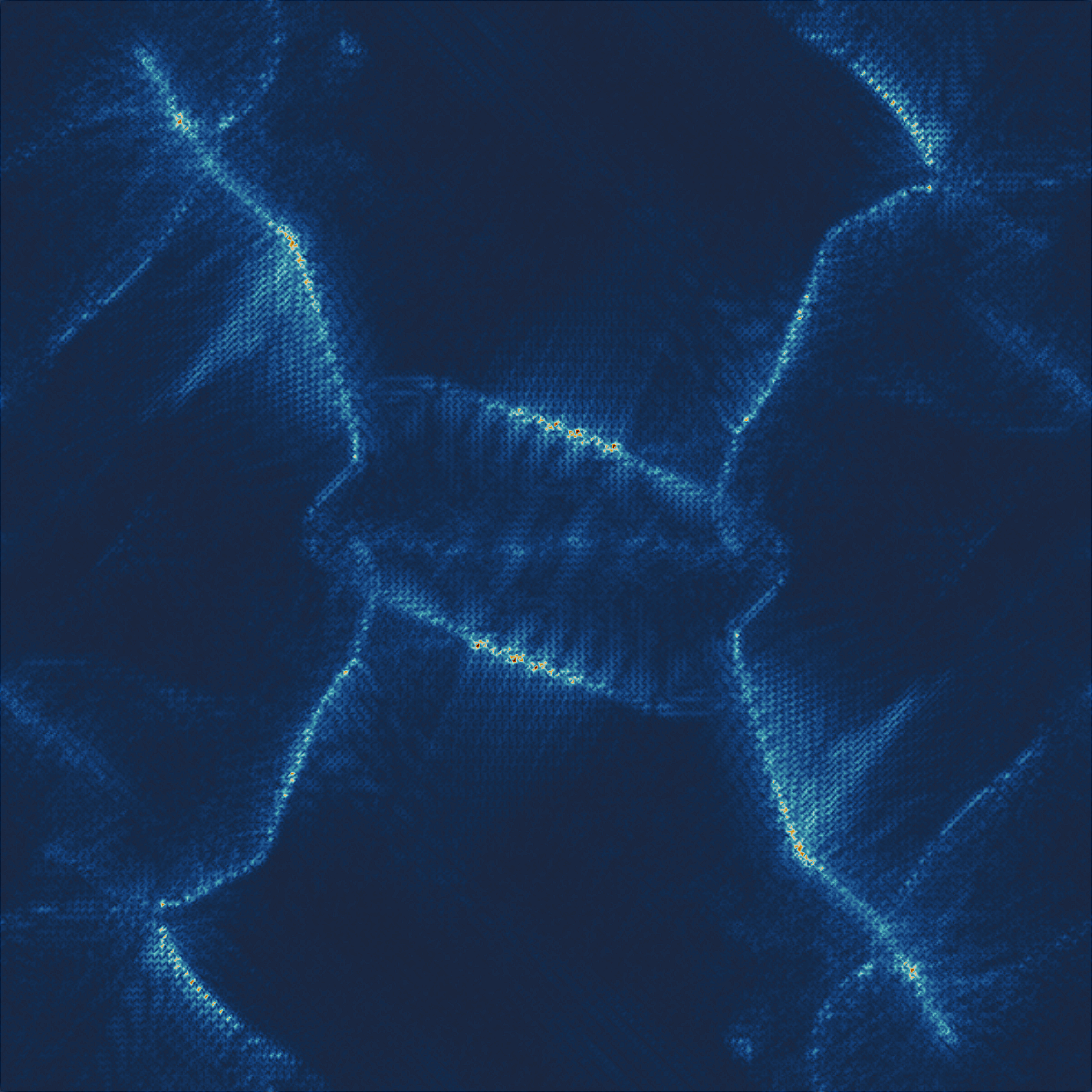}
         \caption{Residual viscosity at $t=0.5$,\\$\mu_h\in[1.71\text{E-}6,3.84\text{E-}4]$}
     \end{subfigure}
     \hfill
     \begin{subfigure}{0.49\textwidth}
         \centering
         \includegraphics[width=\textwidth]{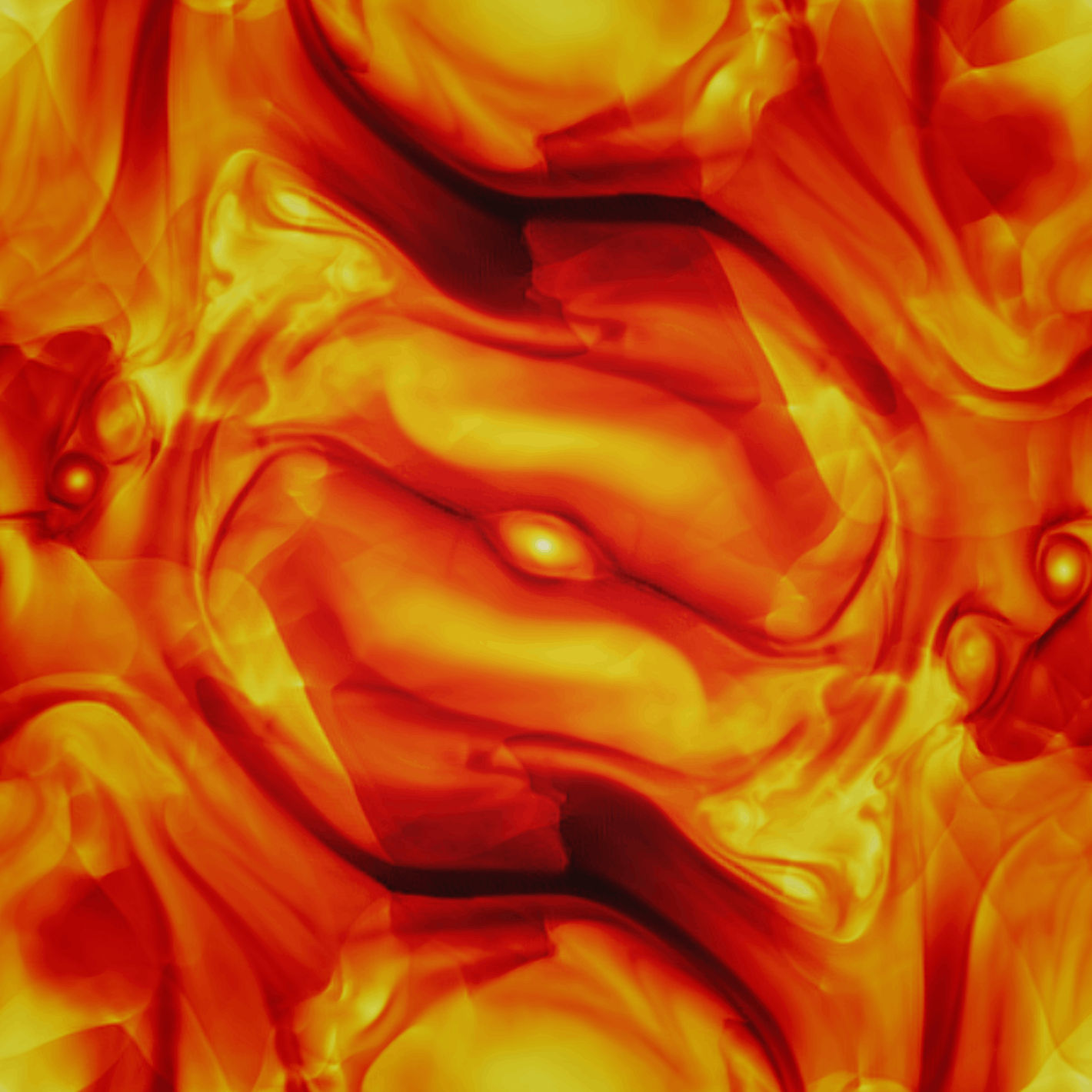}
         \caption{Density at $t=1.0$,\\$\rho_h\in[3.77\text{E-}2,4.19\text{E-}1]$}
     \end{subfigure}
     \hfill
     \begin{subfigure}{0.49\textwidth}
         \centering
         \includegraphics[width=\textwidth]{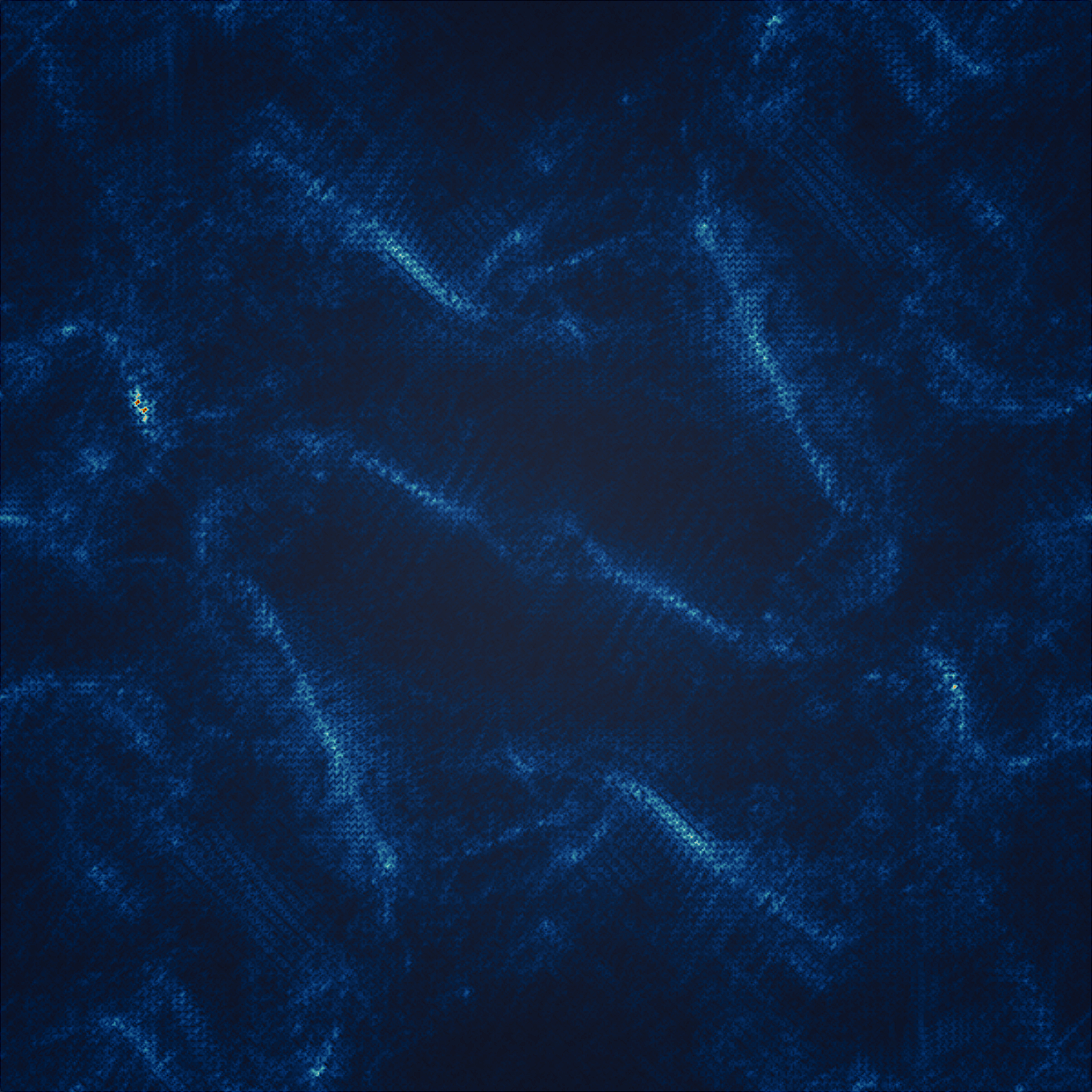}
         \caption{Residual viscosity at $t=1.0$,\\$\mu_h\in[2.67\text{E-}6,3.94\text{E-}4]$}
     \end{subfigure}
     \caption{$\polP_3$ solution of the 2D Orszag-Tang problem, $450\times 450$ nodes}
     \label{fig:OT_2D_P3}
\end{figure}

In this section, we consider the well-known benchmark Orszag-Tang problem in 2D. The considered domain is a square, $\Omega = [0,1] \times [0,1]$. The solution is initialized as
\[
(\rho, \bu, p, \bB)= \left(\frac{25}{36\pi},(-\sin(2\pi y),\sin(2\pi x)),\frac{5}{12\pi},\left(-\frac{\sin(2\pi y)}{\sqrt{4\pi}},\frac{\sin(4\pi x)}{\sqrt{4\pi}}\right)\right).
\]
The adiabatic constant for this problem is $\gamma=\frac{5}{3}$. The periodic boundary condition is used in both directions. For divergence cleaning, we have used the projection method. The density solution $\rho_h$ and the artificial viscosity $\mu_h$ at time $t = 0.5$ and $t=1.0$ are shown in Figure \ref{fig:OT_2D_P1} using $\polP_1$ elements, and in Figure \ref{fig:OT_2D_P3} using $\polP_3$ elements. At $t=1.0$, the structure of the Orszag-Tang solution is considered turbulence \citep{Tricco_2016}.

Despite being a classical benchmark in studies of numerical schemes for MHD, the evolution of the Orszag-Tang solution after $t=0.5$ has rarely been reported or discussed in the literature. For the sake of the readers, we note that the time scale can be different in some other papers. For example, in \citep{Toth_2000}, the solution at $t=3.14$ approximately corresponds to our solution at $t=0.5$. The authors in \citep{Balsara_1998} pointed out that for short integration in time, whether or not a divergence cleaning procedure is involved does not affect the Orszag-Tang solution as there is no noticeable difference. We also observe this situation in our continuous finite element simulation both visually and numerically, at least up to $t=0.5$. In fact, it is evident that the Orszag-Tang solution up to $t=0.5$ is relatively straightforward to be obtained by numerical schemes with shock-capturing capability, e.g., \cite{Balsara_1998, Guillet_2019, Toth_2000}. However, long time integration would discriminate different discretizations, as remarked by \cite{Guillet_2019}. Using DG method, the authors in \cite{Guillet_2019} reported that not all choices of the slope limiters and discretizations of the Powell term \cite{Powell_et_al_1991} would work well after $t = 0.5$. The authors also identified that specifically at $0.75 \leq t \leq 0.85$, numerical solutions are prone to divergence blow-ups due to complex shock collisions.

For this test case, we use the FE setup same to all other test cases without changing any parameter. All the shocks are captured locally and accurately by the RV method. This could be a motivation of using entropy viscosity methods, or in particular RV methods, for shock-capturing purposes. Under the same number of degrees of freedom, the $\polP_3$ solution reveals more structural details of the shocks and turbulence. Similar to the reference \citep{Guillet_2019} using a third-order DG scheme, we capture a twisting density upsurge in the centre of the domain using $\polP_3$ Lagrange elements.

We discuss the effect of different divergence cleaning techniques in the finite element settings in the remaining of this section.


\subsubsection*{Numerical comparison of the divergence cleaning techniques}\label{sec:numerical_div_cleaning_comparison}

\begin{figure}[h!]
     \centering
     \includegraphics[width=0.7\textwidth,viewport=100 520 460 760, clip=true]{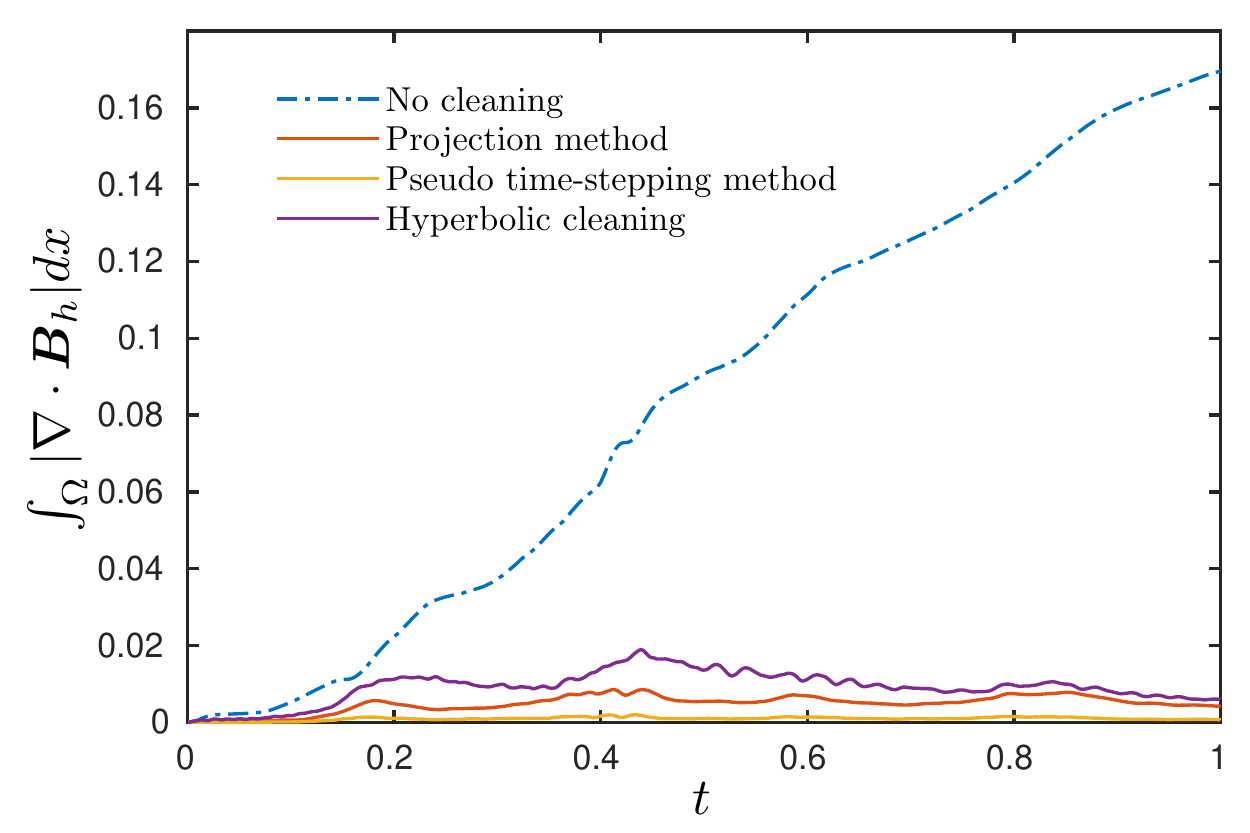}
     \caption{Absolute divergence error $\int_{\Omega}\vert\DIV \bfjlg B_h\vert\ud x$ as a function of time in the 2D Orszag-Tang test. The reference simulation without divergence cleaning corresponds to the blue dot-dashed line. The red, yellow, and purple lines respectively corresponds to the simulations using the projection method, the pseudo time-stepping method, and the hyperbolic cleaning method to clean the divergence.}
     \label{fig:div_clean_compare}
\end{figure}

We show the absolute divergence error $\int_{\Omega}\vert\DIV \bfjlg B_h\vert\ud x$ of the Orszag-Tang test as a function of time using different divergence cleaning techniques in Figure \ref{fig:div_clean_compare}. The simulation uses $\polP_1$ elements on a $40\times 40$ right triangulation mesh. We note that if a sufficiently fine mesh was used, the reference simulation without divergence cleaning would break at some time around $t=0.7$ due to large divergence error. The figure shows that our method works well with some of the most well-known divergence cleaning techniques in MHD. The divergence error remains small with the use of the projection method, pseudo time-stepping method, and hyperbolic method for divergence cleaning. Among the methods, the pseudo time-stepping method with 10 pseudo iterations performs the best. However, the method is the most expensive because 10 pseudo iterations needs to be computed in each time step. The projection method is more expensive than the hyperbolic cleaning method because a Poisson equation needs to be solved in each time step. The obtained behavior of the divergence error agrees with the reported results, e.g., \cite{Tricco_2016, Derigs_2018}. We hereby conclude that the investigated divergence cleaning techniques incorporates well with our shock-capturing method, in which the divergence error remains controlled over long time simulation.

\subsection{3D Orszag-Tang problem, \citep{Helzel_2011}}\label{sec:numerical_OT_3D}
\begin{figure}[h!]
     \centering
     \begin{subfigure}{0.49\textwidth}
         \centering
         \includegraphics[width=\textwidth]{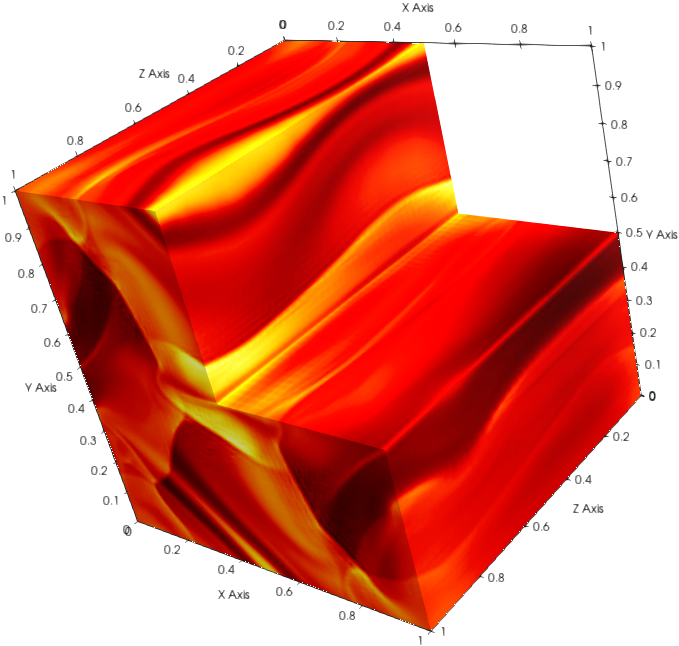}
         \caption{Density at $t=0.5$,\\$\rho_h\in[6.23\text{E-}2,5.79\text{E-}1]$}
     \end{subfigure}
     \hfill
     \begin{subfigure}{0.49\textwidth}
         \centering
         \includegraphics[width=\textwidth]{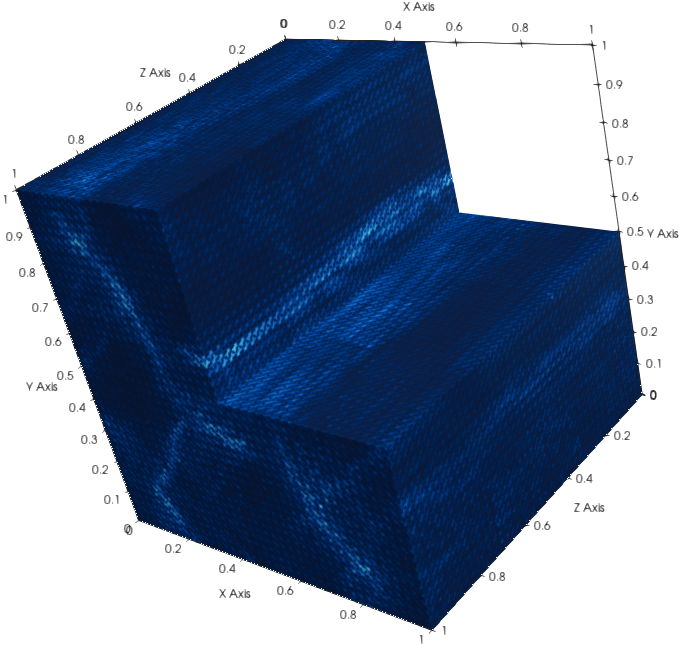}
         \caption{Artificial viscosity at $t=0.5$,\\$\mu_h \in [1.61\text{E-}4,3.65\text{E-}3]$}
     \end{subfigure}
     \caption{$\polP_3$ solution of the 3D Orszag-Tang problem, $120\times 120\times 120$ nodes}
     \label{fig:OT_3D_P3}
\end{figure}

The 3D Orszag-Tang problem is extended from the 2D Orszag-Tang problem in Section \ref{sec:numerical_OT_2D}. Essentially, density and pressure are identical to the 2D problem, independent of the $z$-coordinate. The velocity and magnetic field are also identical in the $x$- and $y$- directions, but zero in the $z$-direction,
\[
(\rho_0, \bu_0, p_0, \bB_0)= \left(\frac{25}{36\pi},(-\sin(2\pi y),\sin(2\pi x),0),\frac{5}{12\pi},\left(-\frac{\sin(2\pi y)}{\sqrt{4\pi}},\frac{\sin(4\pi x)}{\sqrt{4\pi}},0\right)\right).
\]
We then add a small perturbation to the velocity, similar to \cite{Helzel_2011},
\[
\delta \bu = (-\epsilon_p \sin(2\pi z)\sin(2\pi y), \epsilon_p \sin(2\pi z)\sin(2\pi x),\epsilon_p \sin(2\pi z)),
\]
where the perturbation parameter $\epsilon_p$ is a small real number. We set $\epsilon_p = 0.2$ same to \citep{Helzel_2011}.
The initial solution can be written as
\[
(\rho(t), \bu(t), p(t), \bB(t)) =(\rho_0, \bu_0+\delta\bu, p_0, \bB_0+\delta\bB).
\]
The adiabatic constant is kept same to the 2D problem, $\gamma=\frac{5}{3}$. All boundaries are periodic. We use the projection method to clean the divergence. The CFL number is 0.3. A $\polP_3$ solution of the 3D Orszag-Tang problem is shown in Figure~\ref{fig:OT_3D_P3}. We stop the simulation at $t=0.5$ to make the density plot in Figure~\ref{fig:OT_3D_P3}(a) directly comparable with \citep{Helzel_2011}. The density cut $z=0$ presents the highly recognizable 2D Orszag-Tang solution. Along the $z$-axis, the solution agrees well with the reported results in \citep{Helzel_2011}. A quarter of the solution is hidden to reveal some of the inside structures. We see that the RV method can capture the shocks efficiently in 3D with the same finite element setup in 1D and 2D. We omit the discussion of divergence error in 3D as the purpose of this paper is to develop accurate shock-capturing techniques for finite element approximations.

\subsection{MHD Rotor, \cite{Balsara_et_al_1999}} 

\begin{figure}[h!]
     \centering
     \begin{subfigure}{0.49\textwidth}
         \centering
         \includegraphics[width=\textwidth]{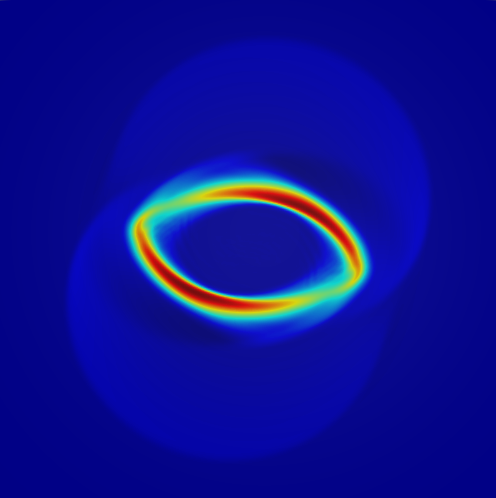}
         \caption{Density at $t=0.15$, $\rho_h\in[7.27\text{E-}1,8.42]$}
     \end{subfigure}
     \hfill
     \begin{subfigure}{0.49\textwidth}
         \centering
         \includegraphics[width=\textwidth]{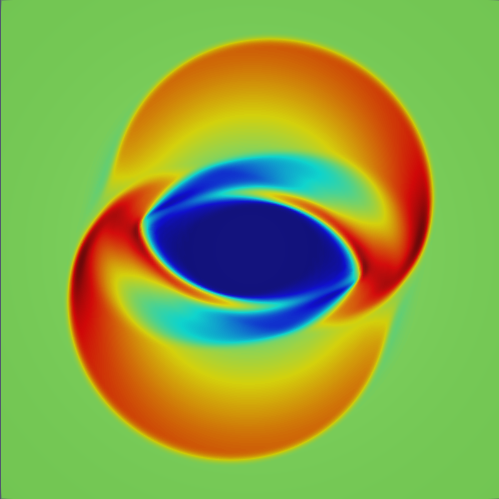}
         \caption{Pressure at $t=0.15$, $p_h\in[3.86\text{E-}2,1.93]$}
     \end{subfigure}
     \hfill
     \begin{subfigure}{0.49\textwidth}
         \centering
         \includegraphics[width=\textwidth]{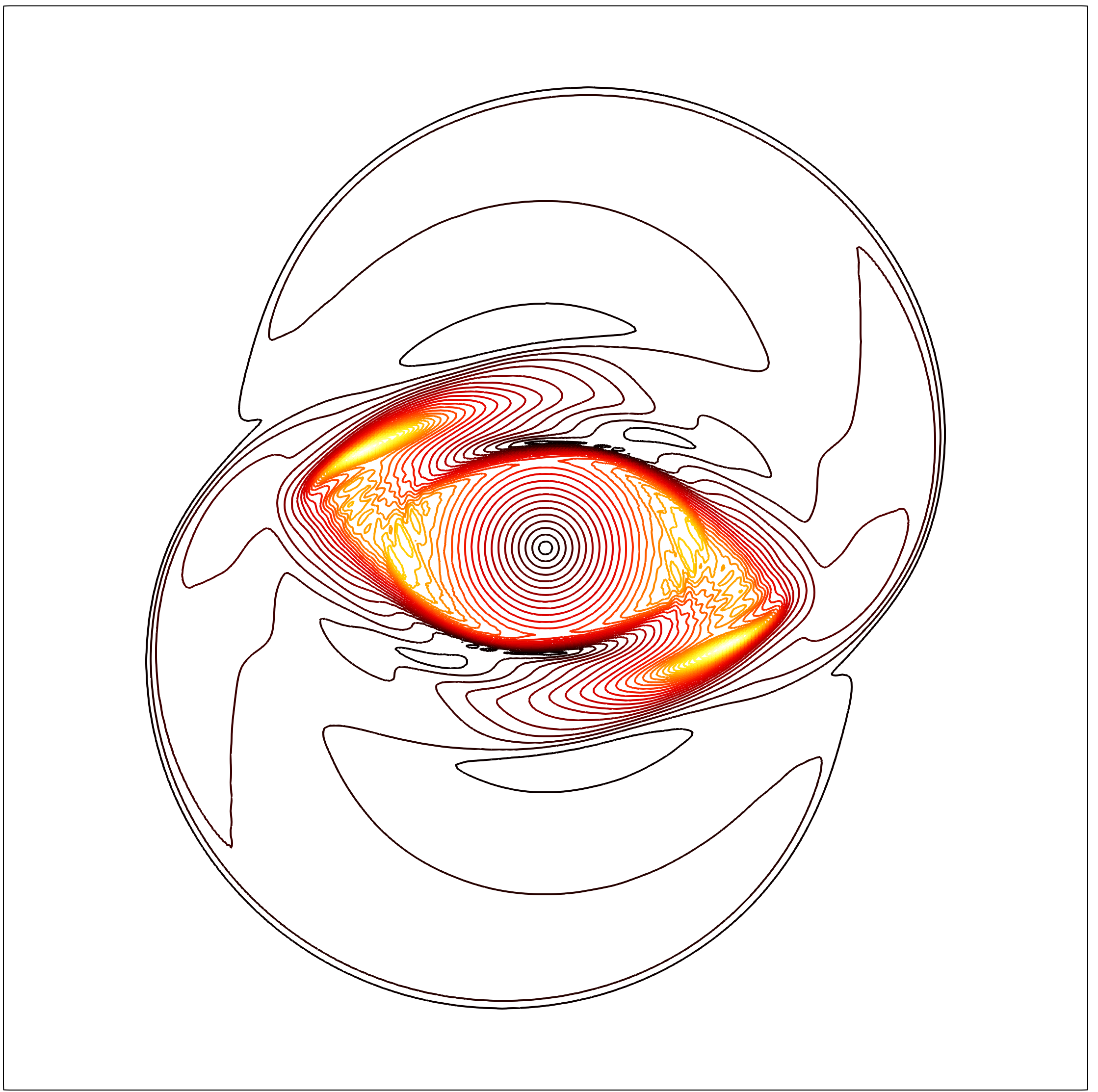}
         \caption{Mach number at $t=0.15$, valued within $[1.67\text{E-}8,4.82]$}
     \end{subfigure}
     \hfill
     \begin{subfigure}{0.49\textwidth}
         \centering
         \includegraphics[width=\textwidth]{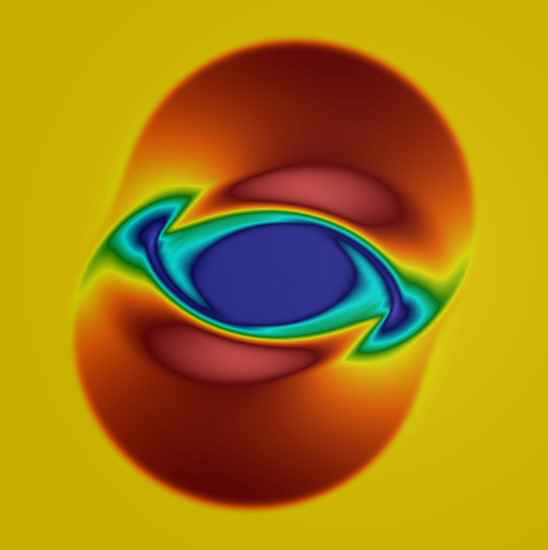}
         \caption{Magnetic pressure at $t=0.15$, $\frac{1}{2}\|\bB_h\|^2\in[5.51\text{E-}2, 2.30]$}
     \end{subfigure}
     \caption{$\polP_3$ solution of the Rotor problem, $300\times 300$ nodes}
     \label{fig:Rotor}
\end{figure}

The final benchmark we present is the rotor problem, first introduced in \cite{Balsara_et_al_1999}. This problem is referred to as the ``first rotor problem'' in \cite{Toth_2000}. The initial solution is a disc being placed in the central of the domain  $\Omega = [0,1] \times [0,1]$. The ambient solution is a stationary medium at ten times lesser density along a homogeneous magnetic field,
\[
(\rho, \bu, p, \bB)= \left(1,(0,0)^T, 1, \left(\frac{5}{\sqrt{4\pi}},0\right)\right).
\]
Let $r$ be the radius of the disc centered at the domain central. The disc is defined as
\[
(\rho, \bu, p)=
\begin{cases}
& \left(10, \left(\frac{u_0}{r_0}(0.5-y),  \frac{u_0}{r_0}(x-0.5)\right)\right) \text{ if }  r < r_0\\
& \left(1+9f, \left(f\frac{u_0}{r}(0.5-y),f\frac{u_0}{r}(x-0.5)\right)\right) \text{ if }r_0 \leq r < r_1,
\end{cases}
\]
where $f = (r_1-r)/(r_1-r_0)$, the inner radius $r_0 = 0.1$ and the outer radius $r_1 = 0.115$. The adiabatic constant is $\gamma=1.4$. Due to imbalance in the centrifugal force, the disc is rotated as it proceeds in time, generating a twisting moment in the magnetic field. The difficulties of this test include: ($i$) resolving the torsional Alfv\'{e}n waves, and ($ii$) during the rotation, the pressure can get close to zero which is challenging for numerical methods. \cite{Toth_2000} reported that many schemes fail to solve this problem owing to negative pressure. Our density and pressure solution at $t=0.15$ are plotted in Figure \ref{fig:Rotor}. The projection method is used for cleaning the divergence. It can be seen that the signature circle contour lines within the middle of Figure~\ref{fig:Rotor}(c) showing the rotational evolution, see e.g., \cite{Toth_2000}, are finely captured.

\section{Conclusion}\label{sec:conclusion}

The main purpose of this paper has been to numerically investigate a new viscous regularization of the ideal MHD equations. The proposed viscous regularization is inspired by the work of \cite{Guermond_2014} on the Euler equations combining with the elliptic terms of the resistive MHD equations. The discretization consists of a high-order residual-based viscosity finite element method in space and the classical explicit Runge-Kuta methods in time. The results have shown that the proposed method can capture different structures of MHD discontinuities while achieving high-order accuracy, experimented up to fourth-order, in shock-free regions.

Future extensions of this paper may concern continuous analysis of the proposed viscous flux regarding e.g., positivity-preserving properties, Galilean invariance, and entropy principles. Discrete invariant-domain preserving schemes can be built upon such continuous properties using, for instance, the viscous operator by \cite{Guermond_Nazarov_2014}.

\backmatter





\bmhead{Acknowledgments}

Some computations were performed on UPPMAX provided by the Swedish National Infrastructure for Computing (SNIC) under project number SNIC 2021/22-233.

\section*{Statements and Declarations}

\bmhead{Funding}
This research is funded by Uppsala University, Sweden.

\bmhead{Competing interests} The authors have no conflicts of interest to disclose.

\bmhead{Data availability statement} Data sharing not applicable to this article as no datasets were generated or analysed during the current study.

\bibliography{ref}
\end{document}